\documentclass{article}
\usepackage[dvipsnames,svgnames,x11names]{xcolor}
\usepackage{arxiv}
\usepackage{amsmath}
\usepackage{amsthm}
\numberwithin{equation}{section}
\usepackage[utf8]{inputenc} % allow utf-8 input
\usepackage[T1]{fontenc}    % use 8-bit T1 fonts
\usepackage{hyperref}       % hyperlinks
\usepackage{url}            % simple URL typesetting
\usepackage{booktabs}       % professional-quality tables
\usepackage{amsfonts}       % blackboard math symbols
\usepackage{nicefrac}       % compact symbols for 1/2, etc.
\usepackage{microtype}      % microtypography
\usepackage{cleveref}       % smart cross-referencing
\usepackage{lipsum}         % Can be removed after putting your text content
\usepackage{graphicx}
\usepackage{natbib}
\usepackage{doi}
\usepackage{tikz}
\usepackage{bm}
\usepackage{indentfirst}
\usepackage{caption}
\usepackage[figuresright]{rotating}
\usepackage{multirow}
\usepackage{threeparttable}
\tikzset{>=stealth}
\usetikzlibrary{patterns}

\title{Parallel generalized solutions of mixed boundary value problem on partially fixed unit annulus subjected to arbitrary traction}

\author{
  \href{https://orcid.org/0000-0000-0000-0000}
  {\includegraphics[scale=0.06]{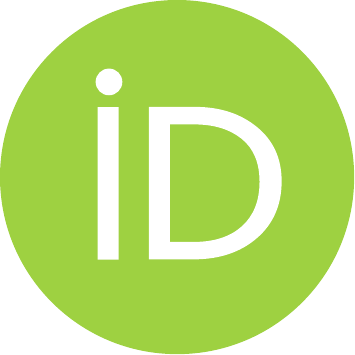}
    \bf {Luobin Lin}} \\
  Fujian Provincial Key Laboratory of Advanced Technology and Informatization in Civil Engineering\\
  College of Civil Engineering\\
  Fujian University of Technology\\
  No. 69 Xueyuan Road, Shangjie University Town, Fuzhou, 350118, Fujian, China \\
  \texttt{luobin\_lin@fjut.edu.cn} \\
  \href{https://orcid.org/0000-0002-5583-3734}
  {\includegraphics[scale=0.06]{orcid.pdf}
    \bf {Fuquan Chen}} \\
  College of Civil Engineering\\
  Fuzhou University\\
  No. 2 Xueyuan Road, Shangjie University Town, Fuzhou, 350108, Fujian, China \\
  \texttt{phdchen@fzu.edu.cn}\\
  \href{https://orcid.org/0000-0003-0679-9127}
  {\includegraphics[scale=0.06]{orcid.pdf}
    \bf {Xianhai Huang}}\\
  College of Civil Engineering\\
  Fujian University of Technology\\
  No. 69 Xueyuan Road, Shangjie University Town, Fuzhou, 350118, Fujian, China \\
  \texttt{hwhsh@163.com} \\
}
\renewcommand{\shorttitle}{}

\begin{document}
\maketitle

\begin{abstract}
  This paper provides two parallel solutions on the mixed boundary value problem of a unit annulus subjected to a partially fixed outer periphery and an arbitrary traction acting along the inner periphery using the complex variable method. The analytic continuation is applied to turn the mixed boundary value problem into a Riemann-Hilbert problem across the free segment along the outer periphery. Two parallel interpreting methods of the unused traction and displacement boundary condition along the outer periphery together with the traction boundary condition along the inner periphery respectively form two parallel complex linear constraint sets, which are then iteratively solved via a successive approximation method to reach the same stable stress and displacement solutions with the Lanczos filtering technique. Finally, four typical numerical cases coded by \texttt{FORTRAN} are carried out and compared to the same cases performed on \texttt{ABAQUS}. The results indicate that these two parallel solutions are both accurate, stable, robust, and fast, and also validate the mutually numerical equivalence of these two parallel solutions.
\end{abstract}

% keywords can be removed
\keywords{Mixed boundary value problem \and Unit annulus \and Riemman-Hilbert problem \and Successive approximation}

\section{Introduction}
\label{sec:intro}

Mixed boundary value problems for elastic annuli are often encountered in composite material, pressure vessel design, oil pipe construction, tunnel engineering, and so on. The merits of linearity makes it feasible that elastic solutions can be presented in analytical or even exact manners to seek profound mechanisms.

The Airy stress function \citep[]{Timoshenko_1951_Elasticity} is a classical method for linearly elastic problems by solving the coefficients of the potential according to boundary conditions, and several analytical solutions on mixed boundary value problems for elastic  lannulus have been provided \cite[]{duffy2008mixed,erdogan1981mixed,belfield1983stress}. Recently, some progress of the Airy stress function method on mixed boundary value problems has been made \cite[]{chawde2021mixed,singh2019simplified} by establishing the simutaneous equations of the strain-displacement relations and the compatibility condtions. 

Compared to the Airy stress function, the complex variable method \cite[]{Muskhelishvili1966, Chenziyin_English} turns to a pair of complex potentials, which are related to the displacement and stress components in close form. By combining with conformal mapping and analytic continuation principle, the complex variable method exhibits powerful ability and flexibility to solve mixed boundary value problems in elastic regions of complicated shapes, especially for simple connected regions \cite[]{hasebe2015mixed, hasebe2021analysis, verma1966mixed, paria1957mixed, ballarini1995certain, hwu1998mixed, selvadurai1985annular, fan1994two, mirsalimov2015crack}.

In this paper, we focus on mixed boundary value problems in an elastic annulus, which is a doubly-connected region. Yau \cite[]{yau1968mixed} proposed a particular solution on a bisymmetrical elastic annulus subjected to fixed constraints acting along two opposite quater arcs of the outer periphery and a constant radial pressure acting along the whole inner periphery. Sugiura \cite[]{Sugiura1973TheMB,sugiura1969heavy} proposed a pair of particular solutions on symmetrical elastic annuli subjected to gravity and axisymmetrically fixed constraints along either periphery. Since then, study on mixed boundary value problem in elastic annulus solved by complex variable method and analytic continuation is rarely seen. The solutions proposed by Yau \cite[]{yau1968mixed} and Sugiura \cite[]{Sugiura1973TheMB,sugiura1969heavy} are elegant, but only consider several simple and symmetrical boundary conditions. Though the results seem correct, ambiguity in the mathematical deduction exists. And due to the limitation of the era, no comparison measures could be used to validate the results.

Therefore, our work consists of the following three parts:

\noindent (a) We propose two parallel generalized solutions on an elastic annulus subjected to a partially fixed constraint acting along the outer periphery and a arbitrary traction acting along the inner periphery by using the complex variable method and analytic continuation principle.

\noindent (b) The mutually numerical equivalence of these two parallel generalized solutions are proven via both analytical deduction and numerical results.

\noindent (c) The analytical results are fully compared with corresponding finite element ones to ensure the validation of the proposed parallel solutions.

\section{Problem definition}
\label{sec:problem-definition}

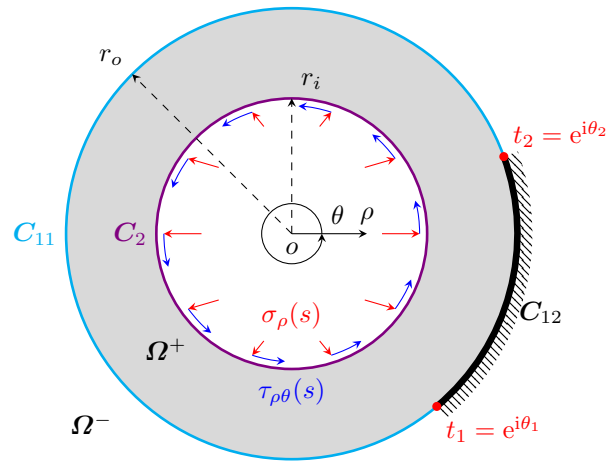
\begin{figure}[htb]
  \centering
  \begin{tikzpicture}
    \filldraw [gray!30] (0,0) circle [radius = 3];
    \draw [line width = 1pt, cyan] (0,0) circle [radius = 3];
    \node at (-3,0) [left, cyan] {$ {\bm C}_{11} $};
    \fill [pattern = north west lines] ({3*cos(20)},{3*sin(20)}) arc [start angle = 20, end angle = -50, radius = 3] -- ({3.2*cos(-50)},{3.2*sin(-50)}) arc [start angle = -50, end angle = 20, radius = 3.2] -- ({3*cos(20)},{3*sin(20)});
    \draw [line width = 2.5pt] ({3*cos(20)},{3*sin(20)}) arc [start angle = 20, end angle = -50, radius = 3];
    \node at ({3*cos(-15)},{3*sin(-15)}) [below right] {$ {\bm C}_{12} $};
    \filldraw [white] (0,0) circle [radius = 1.8];
    \draw [line width = 1pt, violet] (0,0) circle [radius = 1.8];
    \node at (-1.8,0) [left, violet] {$ {\bm C}_{2} $};
    \foreach \x in {0,1,2,...,9} \draw [red, ->] ({1.2*cos(360/10*\x)},{1.5*sin(360/10*\x)}) -- ({1.7*cos(360/10*\x)},{1.7*sin(360/10*\x)});
    \foreach \x in {0,1,2,...,9} \draw [blue,->] ({1.7*cos(360/10*\x)},{1.7*sin(360/10*\x)}) arc [start angle = {360/10*\x}, end angle = {360/10*\x+15}, radius = 1.7];
    \node at ({1.8*cos(-135)},{1.8*sin(-135)}) [below left] {$ {\bm \varOmega}^{+} $};
    \node at ({3.2*cos(-135)},{3.2*sin(-135)}) [below left] {$ {\bm \varOmega}^{-} $};
    \fill [red] ({3*cos(20)},{3*sin(20)}) circle [radius = 0.06];
    \fill [red] ({3*cos(-50)},{3*sin(-50)}) circle [radius = 0.06];
    \node at ({3*cos(-50)},{3*sin(-50)}) [red, below right] {$ t_{1} = {\rm e}^{{\rm i}\theta_{1}} $};
    \node at ({3*cos(20)},{3*sin(20)}) [red, above right] {$ t_{2} = {\rm e}^{{\rm i}\theta_{2}} $};
    \draw (0,0) [->, dashed] -- ({1.8*cos(90)},{1.8*sin(90)}) node [above right] {$ r_{i} $};
    \draw (0,0) [->, dashed] -- ({3*cos(135)},{3*sin(135)}) node [above left] {$ r_{o} $};
    \node at (0,0) [below] {$ o $};
    \draw (0,0) [->] -- (1,0) node [above] {$ \rho $};
    \draw (0.4,0) [->] arc [start angle = 0, end angle = 360, radius = 0.4];
    \node at (0.4,0) [above right] {$ \theta $};
    \node at (0,-1.8) [blue, below] {$ \tau_{\rho\theta}(s) $};
    \node at (0,-1.4) [red, above] {$ \sigma_{\rho}(s) $};
  \end{tikzpicture}
  \caption{Schematic diagram of partially fixed unit annulus subjected to arbitrary traction}
  \label{fig:1}
\end{figure}

Assume a linearly elastic, isotropic, and homogenerous unit annulus is located in the complex plane $ z (\rho o \theta) $, as shown in Fig. \ref{fig:1}. The Poisson's ratio and shear modulus of the annulus are denoted as $ \nu $ and $ G $, respectively. The outer and inner peripheries are denoted as $ {\bm C}_{1} $ and $ {\bm C}_{2} $ with radii of $ r_{o} = 1 $ and $ r_{i} = r $, respectively. The outer periphery $ {\bm C}_{1} $ separates the whole plane into the inner region $ {\bm \varOmega}^{+} ( \rho < r_{o}) $ and the outer region $ {\bm \varOmega}^{-} (\rho > r_{o}) $. The unit annulus takes the region $ r_{i} \leq \rho \leq r_{o} $, which is denoted as $ {\bm \varOmega} $.

According to Muskhelishvili's complex variable method \cite[]{Chenziyin_English,Muskhelishvili1966}, the stress and displacement components in the annulus in polar form can be expressed as:
\begin{subequations}
  \label{eq:2.1}
  \begin{equation}
    \label{eq:2.1a}
    \sigma_{\theta}(z) + \sigma_{\rho}(z) = 2 [ \varphi^{\prime}(z) + \overline{\varphi^{\prime}(z)} ], \quad z = \rho \cdot {\rm e}^{{\rm i}\theta} \in {\bm \varOmega}
  \end{equation}
  \begin{equation}
    \label{eq:2.1b}
    \sigma_{\rho}(z) + {\rm i} \tau_{\rho\theta}(z) = \varphi^{\prime}(z) + \overline{\varphi^{\prime}(z)} - \overline{z} \overline{\varphi^{\prime\prime}(z)} - \frac{\overline{z}}{z}\overline{\psi^{\prime}(z)}, \quad z = \rho \cdot {\rm e}^{{\rm i}\theta} \in {\bm \varOmega}
  \end{equation}
  \begin{equation}
    \label{eq:2.1c}
    g(z) = 2G [u(z) + {\rm i} v(z)] = \kappa\varphi(z) - z \overline{\varphi^{\prime}(z)} - \overline{\psi(z)}, \quad z = \rho \cdot {\rm e}^{{\rm i}\theta} \in {\bm \varOmega}
  \end{equation}
\end{subequations}
where $ \sigma_{\theta} $, $ \sigma_{\rho} $, and $ \tau_{\rho\theta} $ denote hoop, radial, and shear stress components, respectively; $ \varphi(z) $ and $ \psi(z) $ denote complex potentials; the superscripts $ ^{\prime} $ and $ ^{\prime\prime} $ denote first and second deriatives, respectively; the overline above the potentials denotes taking conjugate of corresponding function; $ u $ and $ v $ denote horizontal and vertical displacement components, respectively; $ \kappa $ denotes the Kolosov constant, which is equal to $ 3 - 4 \nu $ and $ (3-\nu)/(1+\nu) $ for plane strain and plane stress conditions, respectively; $ {\rm e} $ and $ {\rm i} $ denote natural logarithmic base and imaginary unit ($ {\rm i}^{2} = -1 $), respectively.

Part of the outer periphery of the annulus is fixed and denoted as $ {\bm C}_{12} $, while the rest is free and denoted as $ {\bm C}_{11} $. The connecting points are marked as $ t_{1} = {\rm e}^{{\rm i}\theta_{1}} $ and $ t_{2} = {\rm e}^{{\rm i}\theta_{2}} $, and $ t_{2} $ is always anti-clockwise to $ t_{1} $. An arbitrary traction is applied along the inner boundary $ {\bm C}_{2} $. Thus, the boundary conditions along the outer and inner peripheries can be written into the following mixed ones as
\begin{subequations}
  \label{eq:2.2}
  \begin{equation}
    \label{eq:2.2a}
    \sigma_{\rho}(t) + {\rm i}\tau_{\rho\theta}(t) = 0, \quad t = {\rm e}^{{\rm i}\theta} \in {\bm C}_{11}
  \end{equation}
  \begin{equation}
    \label{eq:2.2b}
    u(t) + {\rm i}v(t) = 0, \quad t = {\rm e}^{{\rm i}\theta} \in {\bm C}_{12} \\
  \end{equation}
  \begin{equation}
    \label{eq:2.2c}
    \sigma_{\rho}(s) + {\rm i}\tau_{\rho\theta}(s) = f(\sigma), \quad s = r_{i} \cdot {\rm e}^{{\rm i}\theta} \in {\bm C}_{2}
  \end{equation}
\end{subequations}
where
\begin{equation}
  \label{eq:2.3}
   f(\sigma) = f(\theta) = \sum\limits_{k=-\infty}^{\infty} (p_{k}+{\rm i}q_{k}) \cdot {\rm e}^{{\rm i}k\theta}
\end{equation}
$ p_{k} $ and $ q_{k} $ denote radial and tangential coefficients of the traction, which are assumed to be known beforehand. Our problem is to find the particular solution of the complex potentials in Eq. (\ref{eq:2.1}), according to the mixed boundary conditions in Eq. (\ref{eq:2.2}).

\section{Analytic continuation and Riemann-Hilbert problem}
\label{sec:analyt-cont}

To solve the mixed boundary condition problem above, the analytic continuation is adopted to turn the mixed problem into a Riemann-Hilbert problem. To facilitate deduction below, we expand the definition domain of the complex potentials in Eq. (\ref{eq:2.1}) to $ {\bm \varOmega}^{+} $.

Substituting Eq. (\ref{eq:2.1b}) into Eq. (\ref{eq:2.2a}) yields
\begin{equation}
  \label{eq:3.1}
  \varphi^{\prime}(t) = - \overline{\varphi^{\prime}(t)} + \overline{t} \; \overline{\varphi^{\prime\prime}(t)} + \frac{\overline{t}}{t} \; \overline{\psi^{\prime}(t)},\quad t = {\rm e}^{{\rm i}\theta} \in {\bm C}_{11}
\end{equation}
Since the equality $ \overline{t} = \frac{1}{t} $ exists when $ t \in {\bm C}_{11} $, Eq. (\ref{eq:3.1}) turns to
\begin{equation}
  \label{eq:3.2}
  \varphi^{\prime}(t) = - \overline{\varphi^{\prime}}(\frac{1}{t}) + \frac{1}{t} \; \overline{\varphi^{\prime\prime}}(\frac{1}{t}) + \frac{1}{t^{2}} \; \overline{\psi^{\prime}}(\frac{1}{t}),\quad t = {\rm e}^{{\rm i}\theta} \in {\bm C}_{11}
\end{equation}
According to Eq. (\ref{eq:3.2}), the definition domain of $ \varphi^{\prime}(z) $ can be extended from $ {\bm \varOmega}^{+}$ to $ {\bm \varOmega}^{-}$ by crossing the boundary $ {\bm C}_{11} $:
\begin{equation}
  \label{eq:3.3}
  \varphi^{\prime}(z) = - \overline{\varphi^{\prime}}(\frac{1}{z}) + \frac{1}{z} \; \overline{\varphi^{\prime\prime}}(\frac{1}{z}) + \frac{1}{z^{2}} \; \overline{\psi^{\prime}}(\frac{1}{z}),\quad z \in {\bm \varOmega}^{-}
\end{equation}
Since $ \varphi(z) $, $ \psi(z) $, $ \varphi^{\prime} (z) $, and $ \psi^{\prime}(z) $ are all analytic within region $ {\bm \varOmega}^{+} $, and $ \overline{\varphi^{\prime}}(\frac{1}{z}) $, $ \overline{\varphi^{\prime\prime}}(\frac{1}{z}) $, and $ \overline{\psi^{\prime}}(\frac{1}{z}) $ would be subsequently analytic, because $ \frac{1}{z} \in {\bm \varOmega}^{+} $, when $ {z} \in {\bm \varOmega}^{-} $. Thus, the right-hand side of Eq. (\ref{eq:3.3}) would be analytic within the region $ {\bm \varOmega}^{-} $. Subsequently, $ \varphi^{\prime} (\zeta) $ is analytic within the region $ {\bm \varOmega}^{-} $. With Eq. (\ref{eq:3.3}), Eq. (\ref{eq:3.2}) can be rewritten as
\begin{equation}
  \label{eq:3.4}
  \varphi^{\prime+}(t) - \varphi^{\prime-}(t) = 0,\quad t = {\rm e}^{{\rm i}\theta} \in {\bm C}_{11}
\end{equation}
where the superscripts $ ^{+} $ and $ ^{-} $ denote that $ \varphi^{\prime}(z) $ approaches $ {\bm C}_{11} $ from region $ {\bm \varOmega}^{+} $ and $ {\bm \varOmega}^{-} $, respectively.

Replacing $ z $ with $ \overline{z}^{-1} $ in Eq. (\ref{eq:3.3}) and taking conjugate gives
\begin{equation}
  \label{eq:3.5}
  \psi^{\prime} (z) = \frac{1}{z^{2}} \overline{\varphi^{\prime}}(\frac{1}{z}) + \frac{1}{z^{2}} \varphi^{\prime} (z) - \frac{1}{z} \varphi^{\prime\prime}(z), \quad z \in {\bm \varOmega}^{+}
\end{equation}
Eq. (\ref{eq:3.5}) indicates that $ \psi^{\prime}(z) $ is defined and analytic in $ {\bm \varOmega}^{+} $, but is not defined in $ {\bm \varOmega}^{-} $.

Substituting Eq. (\ref{eq:2.1c}) into Eq. (\ref{eq:2.2b}) yields
\begin{equation}
  \label{eq:3.6}
  g(t) = \kappa \varphi(t) - t \overline{\varphi^{\prime} (t)} - \overline{\psi(t)} = 0,\quad t = {\rm e}^{{\rm i}\theta} \in {\bm C}_{12}
\end{equation}
Taking deriative of $ \theta $ on Eq. (\ref{eq:3.6}) ($ t = {\rm e}^{{\rm i}\theta} $) yields
\begin{equation}
  \label{eq:3.7}
  \frac{{\rm d}g(t)}{{\rm d}\theta} = {\rm i}t \kappa \varphi^{\prime} (t) - {\rm i}t \overline{\varphi^{\prime}(t)} + {\rm i}\overline{t} t \overline{\varphi^{\prime\prime}(t)} + {\rm i}\overline{t} \; \overline{\psi^{\prime}(t)} = 0, \quad t = {\rm e}^{{\rm i}\theta} \in {\bm C}_{12}
\end{equation}
Noting that $ \overline{t} = \frac{1}{t} $ along $ {\bm C}_{12} $, and using the chain rule, Eq. (\ref{eq:3.7}) can be written as
\begin{equation}
  \label{eq:3.8}
  \frac{{\rm d}g(t)}{{\rm d}t} = \frac{1}{{\rm i}t} \frac{{\rm d}g(t)}{{\rm d}\theta} = \kappa \varphi^{\prime}(t) - \overline{\varphi^{\prime}}(\frac{1}{t}) + \frac{1}{t} \overline{\varphi^{\prime\prime}}(\frac{1}{t}) + \frac{1}{t^{2}} \overline{\psi^{\prime}}(\frac{1}{t}) = 0, \quad t = {\rm e}^{{\rm i}\theta} \in {\bm C}_{12}
\end{equation}
Apparently, the last three items on the right-hand side of Eq. (\ref{eq:3.8}) are the same to those in Eq. (\ref{eq:3.2}). Considering Eq. (\ref{eq:3.3}), Eq. (\ref{eq:3.8}) can be rewritten as
\begin{equation}
  \label{eq:3.9}
  \kappa \varphi^{\prime+}(t) + \varphi^{\prime-}(t) = 0, \quad t = {\rm e}^{{\rm i}\theta} \in {\bm C}_{12}
\end{equation}
Eqs. (\ref{eq:3.4}) and (\ref{eq:3.9}) form a homogenerous Riemman-Hilbert problem in the whole plane $ 0 < \rho < \infty $.

\section{Solution of the Riemann-Hilbert problem}
\label{sec:solution}

Eq. (\ref{eq:3.9}) suggusts a branch cut along $ {\bm C}_{12} $ in the whole plane, and the following general solution can be found according to Plemelj formula \cite[]{Chenziyin_English,Muskhelishvili1966}:
\begin{equation}
  \label{eq:4.1}
  \varphi^{\prime} (z) = X(z) \sum\limits_{k=-\infty}^{\infty} d_{k}z^{k}, \quad 0 < \rho < \infty
\end{equation}
where 
\begin{equation}
  \label{eq:4.2}
  X(z) = (z - t_{1})^{-\gamma}(z - t_{2})^{\gamma-1}, \quad \gamma = \frac{1}{2} + {\rm i}\lambda, \quad \lambda = \frac{\ln \kappa}{2\pi}, \quad t_{1} = {\rm e}^{{\rm i}\theta_{1}}, t_{2} = {\rm e}^{{\rm i}\theta_{2}}
\end{equation}
$ d_{k} $ are complex coefficients to be determined according to the boundary condition along $ {\bm C}_{2} $ and unused boundary conditions along $ {\bm C}_{1} $. The general solution illustrates poles at origin and infinity, and $ t_{1} $ and $ t_{2} $ are also poles.

Note that the right-hand side of Eq. (\ref{eq:2.3}) is in rational form, thus, $ X(z) $ should be prepared into rational form for coefficient comparisons. $ X(z) $ can be respectively expanded in regions $ {\bm \varOmega}^{+} $ and $ {\bm \varOmega}^{-} $ using Taylor's expansion as
\begin{subequations}
  \label{eq:4.3}
  \begin{equation}
    \label{eq:4.3a}
    X(z) = -t_{1}^{-\gamma} t_{2}^{\gamma-1} (1 - t_{1}^{-1} z)^{-\gamma}(1 - t_{2}^{-1} z)^{\gamma-1} = \sum\limits_{k=0}^{\infty} \alpha_{k} z^{k} , \quad z \in {\bm \varOmega}^{+}
  \end{equation}
  \begin{equation}
    \label{eq:4.3b}
    X(z) = z^{-1} (1 - t_{1}z^{-1})^{-\gamma} (1 - t_{2}z^{-1})^{\gamma - 1} = \sum\limits_{k=1}^{\infty} \beta_{k} z^{-k}, \quad z \in {\bm \varOmega}^{-}
  \end{equation}
\end{subequations}
where the detailed expressions of $ \alpha_{k} $ and $ \beta_{k} $ can be found in Appendix. Substituting Eq. (\ref{eq:4.3}) into Eq. (\ref{eq:4.1}) yields
\begin{subequations}
  \label{eq:4.4}
  \begin{equation}
    \label{eq:4.4a}
    \varphi^{\prime+}(z) = \sum\limits_{k=-\infty}^{\infty} A_{k} z^{k}, \quad z \in {\bm \varOmega}^{+}, \quad A_{k} = \sum\limits_{l=0}^{\infty} \alpha_{l} d_{k-l} 
  \end{equation}
  \begin{equation}
    \label{eq:4.4b}
    \varphi^{\prime-}(z) = \sum\limits_{k=-\infty}^{\infty} B_{k} z^{k}, \quad z \in {\bm \varOmega}^{-}, \quad B_{k} = \sum\limits_{l=1}^{\infty} \beta_{l} d_{k+l} 
  \end{equation}
\end{subequations}
Substituting Eq. (\ref{eq:4.4a}) into Eq. (\ref{eq:3.5}) yields
\begin{equation}
  \label{eq:4.5}
  \psi^{\prime}(z) = \sum\limits_{k=-\infty}^{\infty} \left[ \overline{B}_{-k-2} - (k+1) A_{k+2} \right] z^{k}, \quad  z \in {\bm \varOmega}^{+}
\end{equation}

Substituting Eqs. (\ref{eq:4.4}), (\ref{eq:4.5}) and (\ref{eq:2.1b}) into the inner boundary condition in Eq. (\ref{eq:2.2c}) yields
\begin{equation}
  \label{eq:4.6}
  \sum\limits_{k=-\infty}^{\infty} \left[ A_{k}r^{k} + (k+1)\overline{A}_{-k}(1-r^{-2})r^{-k} - B_{k} r^{-k-2} \right] {\rm e}^{{\rm i}k \theta} = \sum\limits_{k=-\infty}^{\infty} (p_{k}+{\rm i}q_{k}) {\rm e}^{{\rm i}k\theta}
\end{equation}
where $ r $ is the value of the inner radius $ r_{i} $. Comparing the coefficients, we have
\begin{subequations}
  \label{eq:4.7}
  \begin{equation}
    \label{eq:4.7a}
    r^{-1} A_{-1} - r^{-1} B_{-1} = p_{-1}+{\rm i}q_{-1}
  \end{equation}
  \begin{equation}
    \label{eq:4.7b}
    A_{k}r^{k} + (k+1)\overline{A}_{-k}(1-r^{-2})r^{-k} - B_{k} r^{-k-2} = p_{k}+{\rm i}q_{k}, \quad k \leq -2 \; {\rm or} \; k \geq 0
  \end{equation}
\end{subequations}

Since the inner boundary is allowed to deform arbitrarily, the linear constraints on $ d_{k} $ in Eq. (\ref{eq:4.7}) should be complete for the boundary condition in Eq. (\ref{eq:2.2c}). Meanwhile, Eqs. (\ref{eq:3.4}) and (\ref{eq:3.9}) only constrain  the traction and displacement along boundaries $ {\bm C}_{11} $ and $ {\bm C}_{12} $, respectively, but are not completely equivalent to the boundary conditions in Eqs. (\ref{eq:2.2a}) and (\ref{eq:2.2b}). Thus, the displacement and traction along $ {\bm C}_{11} $ and $ {\bm C}_{12}$ should be further examined, respectively, which would provide remaining necessary linear constraints on $ d_{k} $ to form a simultaneous linear constraints.

The displacement along $ {\bm C}_{11} $ and traction along $ {\bm C}_{12}$ can be interpreted via the following two parallel manners, and two parallel generalized solutions are correspondingly presented. Both solutions are satisfactorily accurate for stress and displacement. The major difference between these two solutions lies in the handling of boundary $ {\bm C}_{1} $ and the poles $ t_{1} $ and $ t_{2} $.

\subsection{Solution 1}
\label{sec:solution-1}

Eq. (\ref{eq:3.4}) indicates traction continuation along boundary $ {C}_{11} $, whereas the traction along boundary $ {\bm C}_{12} $ has not been examined. According to static equilibrium for the whole annulus, the resultant acting along the boundary $ {\bm C}_{12} $ should be
\begin{equation}
  \label{eq:4.8}
  \begin{aligned}
    \int_{{\bm C}_{12}} \left[ \sigma_{\rho}(t) + {\rm i}\tau_{\rho\theta}(t) \right] {\rm d}t
    = -\oint_{{\bm C}_{2}} \left[ \sigma_{\rho}(s) + {\rm i}\tau_{\rho\theta}(s) \right] {\rm d}s
    = - {\rm i} \int_{2\pi}^{0} (p_{-1}+{\rm i}q_{-1}) r {\rm d}\theta
    = 2\pi{\rm i}(p_{-1}+{\rm i}q_{-1})r
  \end{aligned}
\end{equation}
The reason that the integral is from $ 2\pi $ to $ 0 $ is to keep the region $ {\bm \varOmega} $ on the left side of the boundary $ {\bm C}_{2} $. On the other hand, 
\begin{equation}
  \label{eq:4.9}
  \begin{aligned}
    \int_{{\bm C}_{12}} \left[ \sigma_{\rho}(t) + {\rm i}\tau_{\rho\theta}(t) \right] {\rm d}t
    & = \int_{{\bm C}_{12}} \left[ \varphi^{\prime+}(t) - \varphi^{\prime-}(t) \right] {\rm d}t \\
    & = \int_{{\bm C}_{12}} \left[ \varphi^{\prime+}(t) - \varphi^{\prime-}(t) \right] {\rm d}t + \int_{{\bm C}_{11}} \left[ \varphi^{\prime+}(t) - \varphi^{\prime-}(t) \right] {\rm d}t \\
    & = \oint_{{\bm C}_{1}} \left[ \varphi^{\prime+}(t) - \varphi^{\prime-}(t) \right] {\rm d}t \\
    & = 2 \pi {\rm i} A_{-1} - 2 \pi {\rm i} B_{-1}
  \end{aligned}
\end{equation}
The last equal sign in Eq. (\ref{eq:4.9}) is due to substitution of Eq. (\ref{eq:4.4}) and residue theorem. Though the definition domains of Eqs. (\ref{eq:4.4a}) and (\ref{eq:4.4b}) are $ z \in {\bm \varOmega}^{+} (0 < \rho <r_{o}) $ and $ z \in {\bm \varOmega}^{-} (\rho > r_{o}) $, respectively, the application premise of residue theorem only considers whether or not the integrand is analytic in a certain region, and the definition domain is not one of the premises. As can be seen in Eqs. (\ref{eq:4.4a}) and (\ref{eq:4.4b}), both $ \varphi^{\prime+}(z) $ and $ \varphi^{\prime-}(z) $ are in rational form, and can be analytic in the whole complex plane, except for the origin and infinity. Therefore, the residue thereom can be applied along $ {\bm C}_{1} $, and the last equal sign in Eq. (\ref{eq:4.9}) stands. Eqs. (\ref{eq:4.8}) and (\ref{eq:4.9}) further analytically validates Eq. (\ref{eq:4.7a}).

Eq. (\ref{eq:3.7}) indicates a constant displacement along boundary $ {\bm C}_{12} $ (zero to be specific), whereas the displacement along boundary $ {\bm C}_{11} $ has not been examined. Apparently, displacements of points $ t_{1} $ and $ t_{2} $ should be equal and zero, thus, we have
\begin{equation}
  \label{eq:4.10}
  \begin{aligned}
    g(t_{1}) - g(t_{2})
    & = \int_{{\bm C}_{11}} [ \kappa \varphi^{\prime+}(t) + \varphi^{\prime-}(t) ] {\rm d}t \\
    & = \int_{{\bm C}_{11}} [ \kappa \varphi^{\prime+}(t) + \varphi^{\prime-}(t) ] {\rm d}t + \int_{{\bm C}_{12}} [ \kappa \varphi^{\prime+}(t) + \varphi^{\prime-}(t) ] {\rm d}t \\
    & = \oint_{{\bm C}_{1}} [ \kappa \varphi^{\prime+}(t) + \varphi^{\prime-}(t) ] {\rm d}t = 0 \\
  \end{aligned}
\end{equation}
Substituting Eq. (\ref{eq:4.4}) into Eq. (\ref{eq:4.10}) with residue theorem yields
\begin{equation}
  \label{eq:4.11}
  \kappa A_{-1} + B_{-1} = 0
\end{equation}
Eq. (\ref{eq:4.11}) indicates single-valueness of displacement in region $ {\bm \varOmega} $.

Solving Eqs. (\ref{eq:4.7a}) and (\ref{eq:4.11}) yields
\begin{equation}
  \label{eq:4.12}
  \begin{cases}
    A_{-1} = \displaystyle \frac {(p_{-1}+{\rm i}q_{-1})r} {1+\kappa} \\
    B_{-1} = \displaystyle \frac {-\kappa (p_{-1}+{\rm i}q_{-1}) r} {1+\kappa} \\
  \end{cases}
\end{equation}
The coefficients in Eq. (\ref{eq:4.12}) coincide with the unbalanced resultant along the inner periphery $ {\bm C}_{2} $. Eqs. (\ref{eq:4.7b}) and (\ref{eq:4.12}) form the simultaneous linear constraints for $ d_{k} $ of Solution 1.

There are many solution methods in theory, and we have tried some of them. For instance, Eqs. (\ref{eq:4.7b}) and (\ref{eq:4.12}) form a complex linear system, and the linear algebra can be applied to reach the direct solution in theory; however, owing to the existence of $ r^{k} $ and $ r^{-k} $, the condition number of the coefficient matrix is generally too large to reach correct solution. Furthermore, nonconvex optimization technique can also be adopted to reach the solution in theory; however, the convergence is not guaranteed due to many possible local minima. After attempts, the successive approximation method \cite[]{sugiura1969heavy} is adopted for its numerical stability and robustness due to a great reduction of the condition number of the coefficient matrices.

Expanding Eqs. (\ref{eq:4.12}) and (\ref{eq:4.7b}) according to Eq. (\ref{eq:4.4}) gives
\begin{subequations}
  \label{eq:4.13}
  \begin{equation}
    \label{eq:4.13a}
    \begin{cases}
      \sum\limits_{l=0}^{\infty} \alpha_{l} d_{-1-l} = \displaystyle \frac {(p_{-1}+{\rm i}q_{-1}) r} {1+\kappa} \\
      \sum\limits_{l=1}^{\infty} \beta_{l} d_{-1+l} = \displaystyle \frac {-\kappa (p_{-1}+{\rm i}q_{-1}) r} {1+\kappa}
    \end{cases}
  \end{equation}
  \begin{equation}
    \label{eq:4.13b}
    \begin{cases}
      r^{-k} \sum\limits_{l=0}^{\infty} \alpha_{l} d_{-k-l} + (-k+1) (1-r^{-2}) r^{k} \sum\limits_{l=0}^{\infty} \overline{\alpha}_{l} \overline{d}_{k-l} - r^{k-2} \sum\limits_{l=1}^{\infty} \beta_{l} d_{-k+l} = p_{-k}+{\rm i}q_{-k}, \quad k \geq 2 \\
      r^{k} \sum\limits_{l=0}^{\infty} \alpha_{l} d_{k-l} + (k+1) (1-r^{-2}) r^{-k} \sum\limits_{l=0}^{\infty} \overline{\alpha}_{l} \overline{d}_{-k-l} - r^{-k-2} \sum\limits_{l=1}^{\infty} \beta_{l} d_{k+l} = p_{k}+{\rm i}q_{k}, \quad k \geq 0
    \end{cases}
  \end{equation}
\end{subequations}
Eq. (\ref{eq:4.13}) can be organized into the following equilibriums:
\begin{subequations}
  \label{eq:4.14}
  \begin{equation}
    \label{eq:4.14a}
    \begin{cases}
      \sum\limits_{l=0}^{\infty} \alpha_{l} d_{-1-l} = \displaystyle \frac {(p_{-1}+{\rm i}q_{-1})r} {1+\kappa} \\
      \sum\limits_{l=0}^{\infty} \alpha_{l} d_{-k-l} = (k-1)(1-r^{-2})r^{2k} \overline{A}_{k} + r^{2k-2} B_{-k} + r^{k}(p_{-k}+{\rm i}q_{-k}), \quad k \geq 2
    \end{cases}
  \end{equation}
  \begin{equation}
    \label{eq:4.14b}
    \begin{cases}
      \sum\limits_{l=1}^{\infty} \beta_{l} d_{-1+l} = & \displaystyle \frac {-\kappa(p_{-1}+{\rm i}q_{-1})r} {1+\kappa} \\
      \sum\limits_{l=1}^{\infty} \beta_{l} d_{l} = & r^{2} A_{0} + (1-r^{-2})r^{2} \overline{A}_{0} + r^{2}(p_{0}+{\rm i}q_{0}) \\
      \sum\limits_{l=1}^{\infty} \beta_{l} d_{1+l} = & r^{4} A_{1} + 2(1-r^{-2})r^{3} \displaystyle \frac {x_{0}-{\rm i}y_{0}} {1+\kappa} + r^{3}(p_{1}+{\rm i}q_{1}) \\
      \sum\limits_{l=1}^{\infty} \beta_{l} d_{k+l} = & r^{2k+2} A_{k} + (k^{2}-1)(1-r^{-2})^{2}r^{2k+2} A_{k} + (k+1)(1-r^{-2})r^{2k} \overline{B}_{-k} \\
                                                      & + (k+1)(1-r^{-2})r^{k+2} (p_{-k}-{\rm i}q_{-k}) + r^{k+2}(p_{k}+{\rm i}q_{k}), \quad k \geq 2
    \end{cases}
  \end{equation}
\end{subequations}
where $ A_{k} (k \geq 0) $ and $ B_{-k} (k \geq 2) $ are computed according to Eq. (\ref{eq:4.4}).

Assume that $ d_{k} $ can be expanded as:
\begin{equation}
  \label{eq:4.15}
  \begin{cases}
    d_{k} = \sum\limits_{q=0}^{\infty} d_{k}^{(q)}, \quad k \geq 0 \\
    d_{-k} = \sum\limits_{q=0}^{\infty} d_{-k}^{(q)}, \quad k \geq 1 \\
  \end{cases}
\end{equation}
We seek approximate solutions for $ d_{k}^{(q)} (k \geq 0) $ and $ d_{-k}^{(q)} (k \geq 1) $ according to Eq. (\ref{eq:4.14}) in the following iterative method. When $q=0$, we set
\begin{subequations}
  \label{eq:4.16}
  \begin{equation}
    \label{eq:4.16a}
    \begin{cases}
      \sum\limits_{l=0}^{\infty} \alpha_{l} d_{-1-l}^{(0)} = \displaystyle \frac {(p_{-1}+{\rm i}q_{-1})r} {1+\kappa} \\
      \sum\limits_{l=0}^{\infty} \alpha_{l} d_{-k-l}^{(0)} = r^{k}(p_{-k}+{\rm i}q_{-k}), \quad k \geq 2
    \end{cases}
  \end{equation}
  \begin{equation}
    \label{eq:4.16b}
    \begin{cases}
      \sum\limits_{l=1}^{\infty} \beta_{l} d_{-1+l}^{(0)} = \displaystyle \frac {-\kappa(p_{-1}+{\rm i}q_{-1})r} {1+\kappa} \\
      \sum\limits_{l=1}^{\infty} \beta_{l} d_{l}^{(0)} = r^{2}(q_{0}+{\rm i}q_{0}) \\
      \sum\limits_{l=1}^{\infty} \beta_{l} d_{1+l}^{(0)} = 2(1-r^{-2})r^{3} \displaystyle \frac {p_{-1}-{\rm i}q_{-1}} {1+\kappa} + r^{3}(p_{1}+{\rm i}q_{1}) \\
      \sum\limits_{l=1}^{\infty} \beta_{l} d_{k+l}^{(0)} = (k+1)(1-r^{-2})r^{k+2} (p_{-k}-{\rm i}q_{-k}) + r^{k+2}(p_{k}+{\rm i}q_{k}), \quad k \geq 2
    \end{cases}
  \end{equation}
\end{subequations}
Eq. (\ref{eq:4.16}) determines $ d_{k}^{(0)} (k \geq 0) $ and $ d_{-k}^{(0)} (k \geq 1) $ to start the iteration. Then for $ q \geq 0 $, $ A_{k}^{(q)} (k \geq 0) $ and $ B_{-k}^{(q)} (k \geq 2) $ can be computed according to Eq. (\ref{eq:4.4}) as:
\begin{equation}
  \label{eq:4.17}
  \begin{cases}
    A_{k}^{(q)} = \sum\limits_{l=0}^{\infty} \alpha_{l} d_{k-l}^{(q)}, \quad k \geq 0 \\
    B_{-k}^{(q)} = \sum\limits_{l=1}^{\infty} \beta_{l} d_{-k+l}^{(q)}, \quad k \geq 2 \\
  \end{cases}
\end{equation}
For $ q \geq 1 $, $ d_{k}^{(q)} (k \geq 0) $ and $ d_{-k}^{(q)} (k \geq 1) $ can be computed as:
\begin{subequations}
  \label{eq:4.18}
  \begin{equation}
  \label{eq:4.18a}
  \begin{cases}
    \sum\limits_{l=0}^{\infty} \alpha_{l} d_{-1-l}^{(q)} = 0 \\
    \sum\limits_{l=0}^{\infty} \alpha_{l} d_{-k-l}^{(q)} = (k-1)(1-r^{-2})r^{2k} \overline{A}_{k}^{(q-1)} + r^{2k-2} B_{-k}^{(q-1)}, \quad k \geq 2
  \end{cases}
\end{equation}
\begin{equation}
  \label{eq:4.18b}
  \begin{cases}
    \sum\limits_{l=1}^{\infty} \beta_{l} d_{-1+l}^{(q)} = 0 \\
    \sum\limits_{l=1}^{\infty} \beta_{l} d_{l}^{(q)} = r^{2} A_{0}^{(q-1)} + (1-r^{-2})r^{2} \overline{A}_{0}^{(q-1)} \\
    \sum\limits_{l=1}^{\infty} \beta_{l} d_{1+l}^{(q)} = r^{4} A_{1}^{(q-1)} \\
    \sum\limits_{l=1}^{\infty} \beta_{l} d_{k+l}^{(q)} = r^{2k+2} A_{k}^{(q-1)} + (k^{2}-1)(1-r^{-2})^{2}r^{2k+2} A_{k}^{(q-1)} + (k+1)(1-r^{-2})r^{2k} \overline{B}_{-k}^{(q-1)}, \quad k \geq 2
  \end{cases}
\end{equation}
\end{subequations}
Then set $ q:q+1 $ into Eq. (\ref{eq:4.17}) to procceed the iteration.

Since $ r \in (0,1) $, all coefficients in the front of $ A_{k}^{(q-1)} (k \geq 0) $ and $ B_{-k}^{(q-1)} (k \geq 2) $ in Eq. (\ref{eq:4.18}) are less than 1, including their conjugates, thus, the right-hand sides of Eqs. (\ref{eq:4.18a}) and (\ref{eq:4.18b}) would gradually approach zero as iteration proceeds. As all elements of the constant vector appoach zero in iteration, the linear solutions would also approach zero. Consequently, the convergence of the iteration procedure in Eqs. (\ref{eq:4.16})-(\ref{eq:4.18}) is guaranteed.

\subsection{Solution 2}
\label{sec:solution-2}

In Eqs. (\ref{eq:4.9}) and (\ref{eq:4.10}) of Solution 1, the application of residue thereom eliminates the two poles $ t_{1} $ and $ t_{2} $ in a mathematically reasonable manner by breaching the definition domains of $ \varphi^{\prime+}(z) $ and $ \varphi^{\prime-}(z) $, and leads a mathematically elegant analytical solution. In Solution 2, we strictly confine $ \varphi^{\prime+}(z) $ and $ \varphi^{\prime-} $ within the definition domains $ {\bm \varOmega}^{+} $ and $ {\bm \varOmega}^{-} $, respectively, without applicationn of residue thereom.

According to Eqs. (\ref{eq:3.9}) and (\ref{eq:3.4}), Eqs. (\ref{eq:4.9}) and (\ref{eq:4.10}) can be respectively rewritten as
\begin{equation}
  \label{eq:4.19}
  \int_{{\bm C}_{12}} \left[ \sigma_{\rho}(t) + {\rm i}\tau_{\rho\theta}(t) \right] {\rm d}t = \int_{{\bm C}_{12}} \left[ \varphi^{\prime+}(t) - \varphi^{\prime-}(t) \right] {\rm d}t = -\frac{1+\kappa}{\kappa} \int_{{\bm C}_{12}} \varphi^{\prime-}(t) {\rm d}t = - \frac{1+\kappa}{\kappa} \int_{{\bm C}_{12}} \varphi^{\prime}(t) {\rm d}t
\end{equation}
\begin{equation}
  \label{eq:4.20}
  g(t_{1}) - g(t_{2}) = \int_{{\bm C}_{11}} \left[ \kappa \varphi^{\prime+}(t) - \varphi^{\prime-}(t) \right] {\rm d}t = (1+\kappa) \int_{{\bm C}_{11}} \varphi^{\prime+}(t) {\rm d}t = (1+\kappa) \int_{{\bm C}_{11}} \varphi^{\prime}(t) {\rm d}t
\end{equation}
In Eq. (\ref{eq:4.19}), $ \varphi^{\prime-}(z) $ is used in integration, instead of $ \varphi^{\prime+}(z) $, because the traction along boundary $ {\bm C}_{12} $ is the constraining force imposed by the region $ {\bm \varOmega}^{-} $. Meanwhile, in Eq. (\ref{eq:4.20}), $ \varphi^{\prime+}(z) $ is used in integration, instead of $ \varphi^{\prime-}(z) $, because the displacement difference of points $ t_{1} $ and $ t_{2} $ should occur within the annulus $ {\bm \varOmega} $, which has an intersection area with region $ {\bm \varOmega}^{+} $. Comparing Eqs. (\ref{eq:4.19}) and (\ref{eq:4.20}) with Eqs. (\ref{eq:4.9}) and (\ref{eq:4.10}), respectively, the definition domains of $ \varphi^{\prime+}(z) $ and $ \varphi^{\prime-}(z) $ are not breached. Such a difference would result in a different solution procedure and require more mathematical labor.

According to static equilibrium of the annulus and Eq. (\ref{eq:4.8}), Eq. (\ref{eq:4.19}) should satisfy
\begin{equation}
  \label{eq:4.21}
  \int_{{\bm C}_{12}} \varphi^{\prime}(t) {\rm d}t = - \frac{2\pi{\rm i}\kappa}{1+\kappa} (p_{-1}+{\rm i}q_{-1})r
\end{equation}
Meanwhile, according to the displacement equality in Eq. (\ref{eq:4.10}), Eq. (\ref{eq:4.20}) should satisfy
\begin{equation}
  \label{eq:4.22}
  \int_{{\bm C}_{11}} \varphi^{\prime}(t) {\rm d}t = 0
\end{equation}
Substituting Eq. (\ref{eq:4.1}) into Eqs. (\ref{eq:4.19}) and (\ref{eq:4.20}) yields
\begin{equation}
  \label{eq:4.23}
  \int_{{\bm C}_{12}} \varphi^{\prime}(t) {\rm d}t = \sum\limits_{k=-\infty}^{\infty} c_{12,k} d_{k}
\end{equation}
\begin{equation}
  \label{eq:4.24}
  \int_{{\bm C}_{11}} \varphi^{\prime}(t) {\rm d}t = \sum\limits_{k=-\infty}^{\infty} c_{11,k} d_{k}
\end{equation}
where
\begin{subequations}
  \label{eq:4.25}
  \begin{equation}
    \label{eq:4.25a}
    c_{12,k} = {\rm i} \int_{\theta_{2}}^{\theta_{1}} ({\rm e}^{{\rm i}\theta}-{\rm e}^{{\rm i}\theta_{1}})^{-\gamma} ({\rm e}^{{\rm i}\theta}-{\rm e}^{{\rm i}\theta_{2}})^{\gamma-1} \cdot {\rm e}^{{\rm i}(k+1)\theta} {\rm d}\theta
  \end{equation}
  \begin{equation}
    \label{eq:4.25b}
    c_{11,k} = {\rm i} \int_{\theta_{2}}^{\theta_{1}+2\pi} ({\rm e}^{{\rm i}\theta}-{\rm e}^{{\rm i}\theta_{1}})^{-\gamma} ({\rm e}^{{\rm i}\theta}-{\rm e}^{{\rm i}\theta_{2}})^{\gamma-1} \cdot {\rm e}^{{\rm i}(k+1)\theta} {\rm d}\theta
  \end{equation}
\end{subequations}
Eq. (\ref{eq:4.25a}) and (\ref{eq:4.25b}) integrates from $ \theta_{2} $ to $ \theta_{1} $ and from $ \theta_{2} $ to $ \theta_{1}+2\pi $ to keep the regions $ {\bm \varOmega}^{-} $ and $ {\bm \varOmega}^{+} $ on the left sides of the boundaries $ {\bm C}_{12} $ and $ {\bm C}_{11} $, respectively. The coefficients in Eq. (\ref{eq:4.25}) can be approximately obtained via numerical integrals, and the details are presented in Appendix \ref{sec:B}.

We would examine the mathematical relationship between Eq. (\ref{eq:4.23}) and (\ref{eq:4.7a}) to see whether it is the same as that in Eq. (\ref{eq:4.9}). Eqs. (\ref{eq:4.7a}) and (\ref{eq:4.23}) can be expanded and rewritten as
\begin{subequations}
  \label{eq:4.26}
  \begin{equation}
    \label{eq:4.26a}
    \sum\limits_{l=0}^{\infty} \alpha_{l}d_{-1-l} + \sum\limits_{l=1}^{\infty} (-\beta_{l})d_{-1+l} = (p_{-1}+{\rm i}q_{-1}) r
  \end{equation}
  \begin{equation}
    \label{eq:4.26b}
    \sum\limits_{l=0}^{\infty} \left(- \frac{1+\kappa}{2\pi{\rm i}\kappa}\right) c_{12,-1-l} d_{-1-l} + \sum\limits_{l=1}^{\infty} \left(- \frac{1+\kappa}{2\pi{\rm i}\kappa}\right) c_{12,-1+l} d_{-1+l} = (p_{-1}+{\rm i}q_{-1}) r
  \end{equation}
\end{subequations}
Comparing the coefficients in Eq. (\ref{eq:4.26}) with consideration of possible multi-valueness of $ \alpha_{l} $ in Eq. (\ref{eq:A.1}) yields
\begin{subequations}
  \label{eq:4.27}
  \begin{equation}
    \label{eq:4.27a}
    c_{12,-1-l} =
    \begin{cases}
      \displaystyle - \frac{2\pi{\rm i}\kappa}{1+\kappa} \alpha_{l}, \quad \theta_{1} \; {\rm and} \; \theta_{2} \; {\rm are \; in \; same \; period} \\
      \displaystyle {\rm e}^{-2\pi\lambda} \cdot \frac{2\pi{\rm i}\kappa}{1+\kappa} \alpha_{l}, \quad \theta_{1} \; {\rm and} \; \theta_{2} \; {\rm are \; in \; different \; periods}\\
    \end{cases}
    , \quad l \geq 0
  \end{equation}
  \begin{equation}
    \label{eq:4.27b}
    c_{12,-1+l} = \displaystyle \frac{2\pi{\rm i}\kappa}{1+\kappa} \beta_{l}, \quad l \geq 1
  \end{equation}
\end{subequations}
In Eq. (\ref{eq:4.27}), one period is $ [-\pi+2n\pi, \pi+2n\pi) $, where $ n $ is an arbitrary integer. The equilities in Eq. (\ref{eq:4.27}) will be numerically examined in the numerical cases in Section \ref{sec:numerical-cases}. Thus, Eq. (\ref{eq:4.23}) is equivalent to Eq. (\ref{eq:4.7a}), which is exactly same as the result derived from Eq. (\ref{eq:4.9}).

Eqs. (\ref{eq:4.24}) and (\ref{eq:4.7}) make up the simultaneous linear constraint system for $ d_{k} $ of Solution 2, which is verified to be equivalent to the simultaneous linear system of Solution 1 by the numerical cases in Section \ref{sec:numerical-cases}. Similar to Solution 1, the successive approximation method is applied. To ensure convergence, we slightly alter the approximation strategy of Solution 1. Eq. (\ref{eq:4.7b}) can be expanded and rewritten as
\begin{subequations}
  \label{eq:4.28}
  \begin{equation}
    \label{eq:4.28a}
    \sum\limits_{l=0}^{\infty} \alpha_{l} d_{-k-l} = (k-1)(1-r^{-2})r^{2k} \overline{A}_{k} + r^{2k-2} B_{-k} + r^{k}(p_{-k}+{\rm i}q_{-k}), \quad k \geq 2
  \end{equation}
  \begin{equation}
    \label{eq:4.28b}
    \begin{cases}
      \sum\limits_{l=1}^{\infty} \beta_{l} d_{l} = & r^{2} A_{0} + (1-r^{-2})r^{2} \overline{A}_{0} + r^{2}(p_{0}+{\rm i}q_{0}) \\
      \sum\limits_{l=1}^{\infty} \beta_{l} d_{1+l} = & r^{4} A_{1} + 2(1-r^{-2})r^{2} \overline{B}_{-1} + 2(1-r^{-2})r^{3}(p_{-1}-{\rm i}q_{-1}) + r^{3}(p_{1}+{\rm i}q_{1}) \\
      \sum\limits_{l=1}^{\infty} \beta_{l} d_{k+l} = & r^{2k+2} A_{k} + (k^{2}-1)(1-r^{-2})^{2}r^{2k+2} A_{k} + (k+1)(1-r^{-2})r^{2k} \overline{B}_{-k} \\
                                                   & + (k+1)(1-r^{-2})r^{k+2}(p_{-k}-{\rm i}q_{-k}) + r^{k+2}(p_{k}+{\rm i}q_{k}), \quad k \geq 2
    \end{cases}
  \end{equation} 
\end{subequations}
While Eq. (\ref{eq:4.7a}) and (\ref{eq:4.25}) can be expanded and rewritten as
\begin{subequations}
  \label{eq:4.29}
  \begin{equation}
    \label{eq:29a}
    \alpha_{0} d_{-1} - \beta_{1} d_{0} = - \sum\limits_{l=1}^{\infty} \alpha_{l} d_{-1-l} + \sum\limits_{l=2}^{\infty} \beta_{l} d_{-1+l} + r (p_{-1} + {\rm i}q_{-1})
  \end{equation}
  \begin{equation}
    \label{eq:29b}
    c_{11,-1} d_{-1} + c_{11,0} d_{0} = - \sum\limits_{l=2}^{\infty} c_{11,-l} d_{-l} - \sum\limits_{l=1}^{\infty} c_{11,l} d_{l}
  \end{equation}
\end{subequations}

Assume that $ d_{k} $ can be expanded into the same form in Eq. (\ref{eq:4.16}), then the following iterative method is applied. When $ q = 0 $, we set
\begin{subequations}
  \label{eq:4.30}
  \begin{equation}
    \label{eq:4.30a}
    \sum\limits_{l=0}^{\infty} \alpha_{l} d_{-k-l}^{(0)} = r^{k}(p_{-k}+{\rm i}q_{-k}), \quad k \geq 2
  \end{equation}
  \begin{equation}
    \label{eq:4.30b}
    \begin{cases}
      \sum\limits_{l=1}^{\infty} \beta_{l} d_{l}^{(0)} = r^{2}(p_{0}+{\rm i}q_{0}) \\
      \sum\limits_{l=1}^{\infty} \beta_{l} d_{k+l}^{(0)} = (k+1)(1-r^{-2})r^{k+2}(p_{-k}-{\rm i}q_{-k}) + r^{k+2}(p_{k}+{\rm i}q_{k}), \quad k \geq 1
    \end{cases}
  \end{equation}
\end{subequations}
Eq. (\ref{eq:4.30}) gives $ d_{-k}^{(0)} (k \geq 2) $ and $ d_{k}^{(0)} (k \geq 1) $, and we can compute $ d_{-1}^{(0)} $ and $ d_{0}^{(0)} $ via
\begin{subequations}
  \label{eq:4.31}
  \begin{equation}
    \label{eq:31a}
    \alpha_{0} d_{-1}^{(0)} - \beta_{1} d_{0}^{(0)} = - \sum\limits_{l=1}^{\infty} \alpha_{l} d_{-1-l}^{(0)} + \sum\limits_{l=2}^{\infty} \beta_{l} d_{-1+l}^{(0)} + r (p_{-1} + {\rm i}q_{-1})
  \end{equation}
  \begin{equation}
    \label{eq:4.31b}
    c_{11,-1} d_{-1}^{(0)} + c_{11,0} d_{0}^{(0)} = - \sum\limits_{l=2}^{\infty} c_{11,-l} d_{-l}^{(0)} - \sum\limits_{l=1}^{\infty} c_{11,l} d_{l}^{(0)}
  \end{equation}
\end{subequations}
Eqs. (\ref{eq:4.30}) and (\ref{eq:4.31}) sequentially determine $ d_{-k}^{(0)} (k \geq 2) $ and $ d_{k}^{(0)} (k \geq 1) $, $ d_{-1}^{(0)} $ and $ d_{0}^{(0)} $ to start iteration. For $ q \geq 0 $, $ A_{k}^{q} (k \ge 0) $ and $ B_{-k}^{(q)} (k \geq 1) $ can be computed as
\begin{equation}
  \label{eq:4.32}
  \begin{cases}
    A_{k}^{(q)} = \sum\limits_{l=0}^{\infty} \alpha_{l} d_{k-l}^{(q)}, \quad k \geq 0 \\
    B_{-k}^{(q)} = \sum\limits_{l=1}^{\infty} \beta_{l} d_{-k+l}^{(q)}, \quad k \geq 1 \\
  \end{cases}
\end{equation}
For $ q \geq 1 $, $ d_{-k}^{(q)} (k \geq 2) $ and $ d_{k}^{(q)} (k \geq 1) $ can be determined as
\begin{subequations}
  \label{eq:4.33}
  \begin{equation}
    \label{eq:4.33a}
    \sum\limits_{l=0}^{\infty} \alpha_{l} d_{-k-l}^{(q)} = (k-1)(1-r^{-2})r^{2k} \overline{A}_{k}^{(q-1)} + r^{2k-2} B_{-k}^{(q-1)}, \quad k \geq 2
  \end{equation}
  \begin{equation}
    \label{eq:4.33b}
    \begin{cases}
      \sum\limits_{l=1}^{\infty} \beta_{l} d_{l}^{(q)} = & r^{2} A_{0}^{(q-1)} + (1-r^{-2})r^{2} \overline{A}_{0}^{(q-1)} \\
      \sum\limits_{l=1}^{\infty} \beta_{l} d_{1+l}^{(q)} = & r^{4} A_{1}^{(q-1)} + 2(1-r^{-2})r^{2} \overline{B}_{-1}^{(q-1)} \\
      \sum\limits_{l=1}^{\infty} \beta_{l} d_{k+l}^{(q)} = & r^{2k+2} A_{k}^{(q-1)} + (k^{2}-1)(1-r^{-2})^{2}r^{2k+2} A_{k}^{(q-1)} + (k+1)(1-r^{-2})r^{2k} \overline{B}_{-k}^{(q-1)}, \quad k \geq 2
    \end{cases}
  \end{equation}
\end{subequations}
Then we can compute $ d_{-1}^{(q)} $ and $ d_{0}^{(q)} $ via
\begin{subequations}
  \label{eq:4.34}
  \begin{equation}
    \label{eq:4.34a}
    \alpha_{0} d_{-1}^{(q)} - \beta_{1} d_{0}^{(q)} = - \sum\limits_{l=1}^{\infty} \alpha_{l} d_{-1-l}^{(q)} + \sum\limits_{l=2}^{\infty} \beta_{l} d_{-1+l}^{(q)}
  \end{equation}
  \begin{equation}
    \label{eq:4.34b}
    c_{11,-1} d_{-1}^{(q)} + c_{11,0} d_{0}^{(q)} = - \sum\limits_{l=2}^{\infty} c_{11,-l} d_{-l}^{(q)} - \sum\limits_{l=1}^{\infty} c_{11,l} d_{l}^{(q)}
  \end{equation}
\end{subequations}
Then set $ q:q+1 $ into Eq. (\ref{eq:4.31}) to proceed the iteration.

Similar to Solution 1, all the coefficients in front of $ A_{k} (k \geq 0) $ and $ B_{-k} (k \geq 1) $ are less than 1, including their conjugates, thus, the right-hand sides of Eq. (\ref{eq:4.33}) would approach zero as iteration proceeds, as well as $ d_{-k}^{(q)} (k \geq 2) $ and $ d_{k}^{(q)} (k \geq 1) $ on the left-hand sides. Subsequently, $ d_{-1}^{(q)} $ and $ d_{0}^{(q)} $ would approach zero, as iteration proceeds. Consequently, the convergence of the iteration procedure in Eq. (\ref{eq:4.30})-(\ref{eq:4.34}) is guaranteed.

\subsection{Final solution with truncation}
\label{sec:numer-solut-final}

To obtain actual computation results, we have to truncate the infinite series in Eqs. (\ref{eq:2.3}) and (\ref{eq:4.1}) into $ 2N+1 $ items. For Solution 1, Eqs. (\ref{eq:4.16a}) and (\ref{eq:4.18a}) turn to simultaneous complex linear systems containing $ N $ complex variables and $ N $ complex linear equations; Eqs. (\ref{eq:4.16b}) and (\ref{eq:4.18b}) turn to simultaneous complex linear systems containing $ N+1 $ complex variables and $N+1 $ complex linear equations; Eq. (\ref{eq:4.18}) becomes finite as well. For Solution 2, Eqs. (\ref{eq:4.30a}) and (\ref{eq:4.33a}) turn to simultaneous complex linear systems containing $ N-1 $ complex variables and $ N-1 $ complex linear equations; Eqs. (\ref{eq:4.30b}) and (\ref{eq:4.33b}) turn to simultaneous complex linear systems containing $ N $ complex variables and $ N $ complex linear equations; Eq. (\ref{eq:4.32}) becomes finite as well. The coefficient matrices in Eqs. (\ref{eq:4.16}) and (\ref{eq:4.18}), (\ref{eq:4.30}) and (\ref{eq:4.33}) are respectively the same and would not be altered in iterations, and the condition number is small as illustrated in the numerical cases.

The iteration may stop when
\begin{equation}
  \label{eq:4.35}
  \max | d_{k}^{(q)} | \leq \epsilon, \quad -N \leq k \leq N
\end{equation}
where $ \epsilon $ denotes the error tolerance. When the iteration stops, the maximum iteration rep is recorded as $ Q $.

The solution $ d_{k} $ gives $ A_{k} $ and $ B_{k} $ in Eq. (\ref{eq:4.4}), and the complex potentials within the annulus can be obtained via Eqs. (\ref{eq:4.4a}) and (\ref{eq:4.5}) as
\begin{subequations}
  \label{eq:4.36}
  \begin{equation}
    \label{eq:4.36a}
    \varphi^{\prime}(z) = \sum\limits_{k=-N}^{N} A_{k} z^{k} \cdot F(k), \quad z \in {\bm \varOmega}, \quad A_{k} = \sum\limits_{l=0}^{N+k} \alpha_{l} d_{k-l} 
  \end{equation}
  \begin{equation}
    \label{eq:4.36b}
    \psi^{\prime}(z) = \sum\limits_{k=-N}^{N} \left[ \overline{B}_{-k} - (k-1) A_{k} \right] z^{k-2} \cdot F(k), \quad z \in {\bm \varOmega}, \quad B_{k} = \sum\limits_{l=1}^{N-k} \beta_{l} d_{k+l} 
  \end{equation}
\end{subequations}
where $ F(k) (-N \leq k \leq N) $ denote the Lanczos filtering parameters \cite[]{Lanczos1956,singh2019simplified,chawde2021mixed}, and can be expressed as
\begin{equation}
  \label{eq:4.37}
  F(k) =
  \begin{cases}
    1, \quad k=0 \\
    \sin\left(\frac{|k|}{N} \pi \right)/(\frac{|k|}{N} \pi), \quad {\rm otherwise} \\
  \end{cases}
\end{equation}
Finally, the stress and displacement components in Eq. (\ref{eq:2.1}) after normalization can be expressed as
\begin{subequations}
  \label{eq:4.38}
  \begin{equation}
    \label{eq:4.38a}
    \Sigma_{\theta}(\rho,\theta) + \Sigma_{\rho}(\rho,\theta) = 4 \Re \sum\limits_{k=-N}^{N} \frac{F(k)}{G} \cdot A_{k} \rho^{k} {\rm e}^{{\rm i}k\theta}, \quad z \in {\bm \varOmega}
  \end{equation}
  \begin{equation}
    \label{eq:4.38b}
    \Sigma_{\rho}(\rho,\theta) + {\rm i} \Sigma_{\rho\theta}(\rho,\theta) = \sum\limits_{k=-N}^{N} \frac{F(k)}{G} \cdot \left[ \rho^{k} A_{k} + (k+1)(1-\rho^{-2})\rho^{-k} \overline{A}_{-k} - \rho^{-k-2} B_{k} \right] {\rm e}^{{\rm i}k \theta}, \quad z \in {\bm \varOmega}
  \end{equation}
  \begin{equation}
    \label{eq:4.38c}
    \begin{aligned}
      U(\rho,\theta)+{\rm i}V(\rho,\theta) = & \sum\limits_{k=0}^{N} \frac{F(k)}{2Gr_{o}} \cdot \left[ \kappa A_{k} \frac{\rho^{k+1}}{k+1} - \overline{A}_{-k}\rho^{-k+1}(1-\rho^{-2}) - B_{k}\frac{\rho^{-k-1}}{-k-1} \right] {\rm e}^{{\rm i}(k+1)\theta} \\
                       & + \frac{F(1)}{2Gr_{o}} \cdot [ ( \kappa A_{-1} - B_{-1} )\log\rho-\overline{A}_{1}\rho^{2}] - (D_{x}+{\rm i}D_{y}) \\
                       & + \sum\limits_{k=2}^{N} \frac{F(k)}{2Gr_{o}} \cdot \left[ \kappa A_{-k} \frac{\rho^{-k+1}}{-k+1} -\overline{A}_{k}\rho^{k+1}(1-\rho^{-2}) - B_{-k}\frac{\rho^{k-1}}{k-1} \right] {\rm e}^{{\rm i}(-k+1)\theta} \\
    \end{aligned}
    , \quad z \in {\bm \varOmega}
  \end{equation}
\end{subequations}
where $ \Sigma_{\theta} = \sigma_{\theta}/G $, $ \Sigma_{\rho} = \sigma_{\rho}/G $, $ \Sigma_{\rho\theta} = \tau_{\rho\theta}/G $, $ U = u/r_{o} $, $ V = v/r_{o} $, $ D_{x} $ and $ D_{y} $ denote horizontal and vertical rigid-body displacements after normalization, respectively.

\section{Numerical cases}
\label{sec:numerical-cases}

We will examine these two parallel solutions above via four numerical cases of unit annuli in plane strain condition. The solutions in these four cases are coded by \texttt{FORTRAN}, and performed on \texttt{GCC 11}. The condition numbers of the coefficient matrice in Eq. (\ref{eq:4.17}) and (\ref{eq:4.19}) are computed using \texttt{ZGESVD} package of \texttt{LAPAPCK/complex16}, and the corresponding complex linear systems are solved using \texttt{ZGESV} package. The error tolerance in Eq. (\ref{eq:4.35}) takes $ \epsilon = 10^{-20} $. For comparisons, the same cases are conducted in \texttt{ABAQUS 2020} for computation using finite element method.

\begin{figure}[htb]
  \centering
  \begin{tabular}{cc}
    \begin{tikzpicture}
      \filldraw [gray!30] (0,0) circle [radius = 3];
      \draw [line width = 1pt, cyan] (0,0) circle [radius = 3];
      \fill [pattern = north west lines] ({3*cos(0)},{3*sin(0)}) arc [start angle = 0, end angle = -90, radius = 3] -- ({3.2*cos(-90)},{3.2*sin(-90)}) arc [start angle = -90, end angle = 0, radius = 3.2] -- ({3*cos(0)},{3*sin(0)});
      \draw [line width = 2.5pt] ({3*cos(0)},{3*sin(0)}) arc [start angle = 0, end angle = -90, radius = 3];
      \filldraw [white] (0,0) circle [radius = 0.3];
      \draw [line width = 1pt, violet] (0,0) circle [radius = 0.3];
      \foreach \x in {0,1,2,...,9} \draw [->] ({0.3*cos(360/10*\x)},{0.3*sin(360/10*\x)}) -- ({0.3*cos(360/10*\x)},{0.3*sin(360/10*\x)+0.3});
      \fill [red] ({3*cos(0)},{3*sin(0)}) circle [radius = 0.06];
      \fill [red] ({3*cos(-90)},{3*sin(-90)}) circle [radius = 0.06];
      \node at ({3*cos(-90)},{3*sin(-90)}) [red, above] {$ \theta_{1} = -\frac{\pi}{2} $};
      \node at ({3*cos(0)},{3*sin(0)}) [red, left] {$ \theta_{2} = 0 $};
      \node at (0,-0.3) [below, align=center] {$ f(\theta) = {\rm i}\cos\theta + \sin\theta $};
      \node at (0,-3.2) [below] {$ r_{i}/r_{o} $ = 0.1, $r_{\rho}/r_{o}$ = 0.3};
    \end{tikzpicture}
    &
      \includegraphics{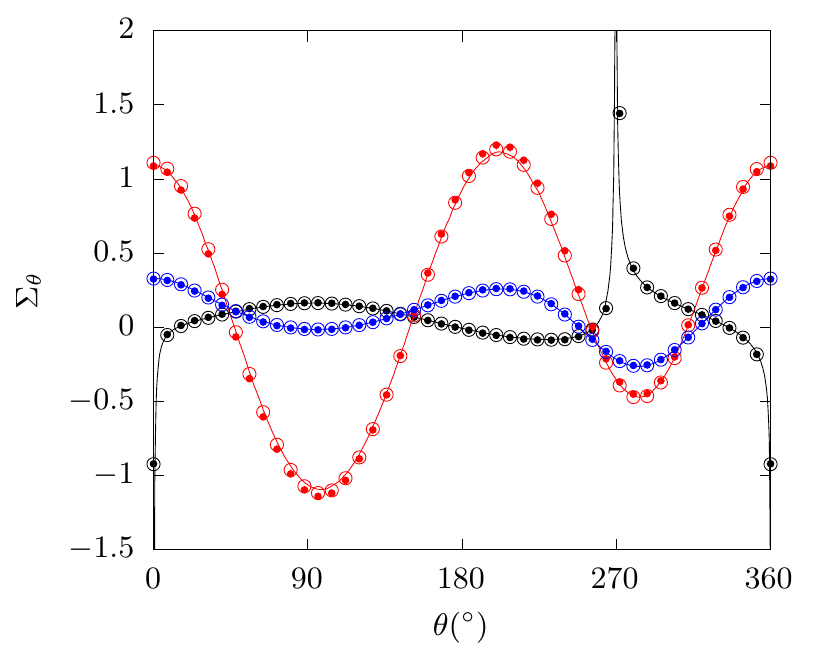}
    \\
    (a) Schematic diagram of Case A
    &
      (b) Hoop stress
    \\
    \includegraphics{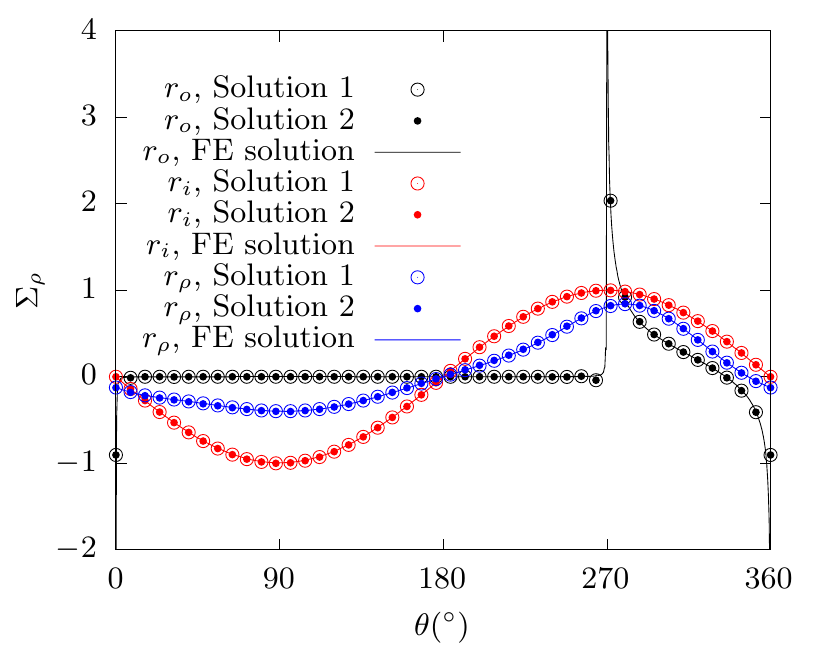}
    &
      \includegraphics{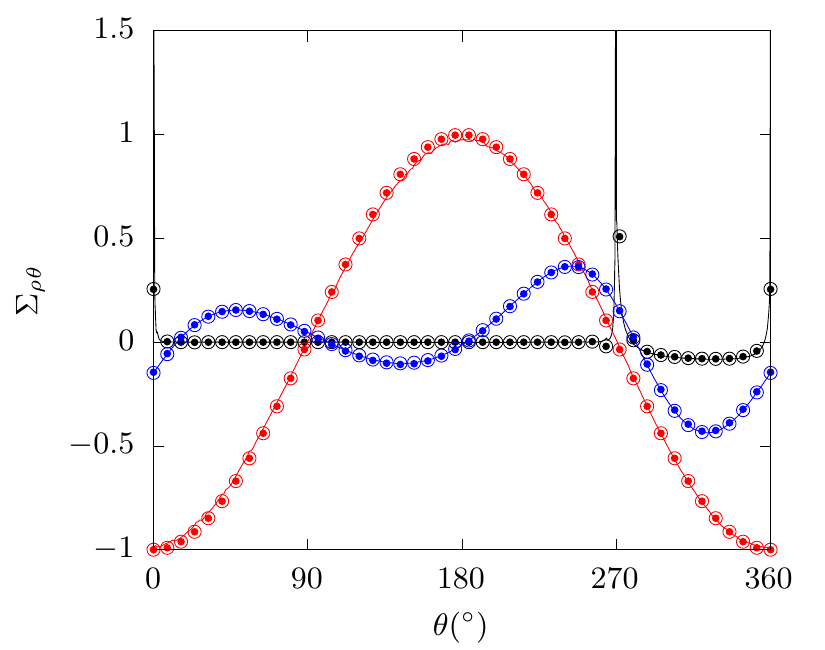}
    \\
    (c) Radial stress
    &
      (d) Shear stress
    \\
    \includegraphics{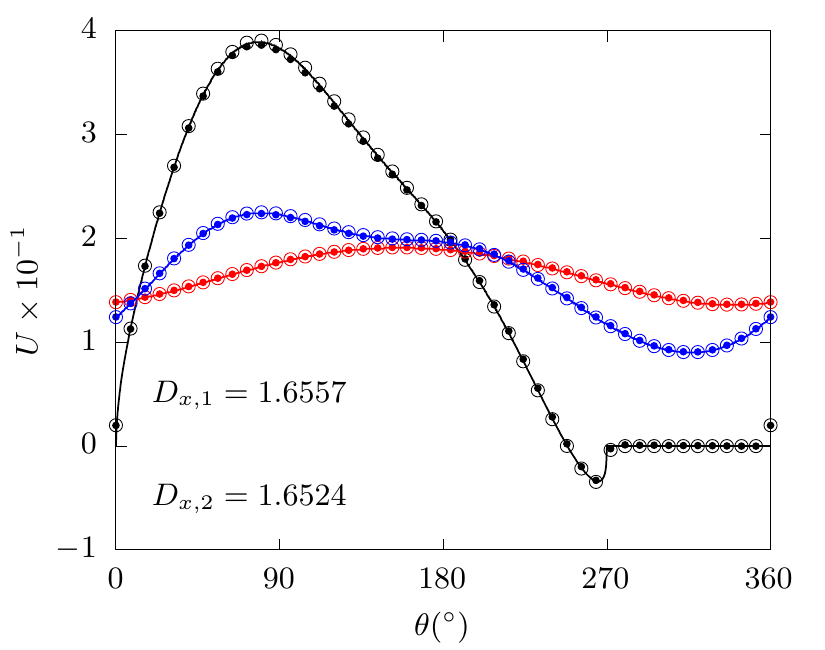}
    &
      \includegraphics{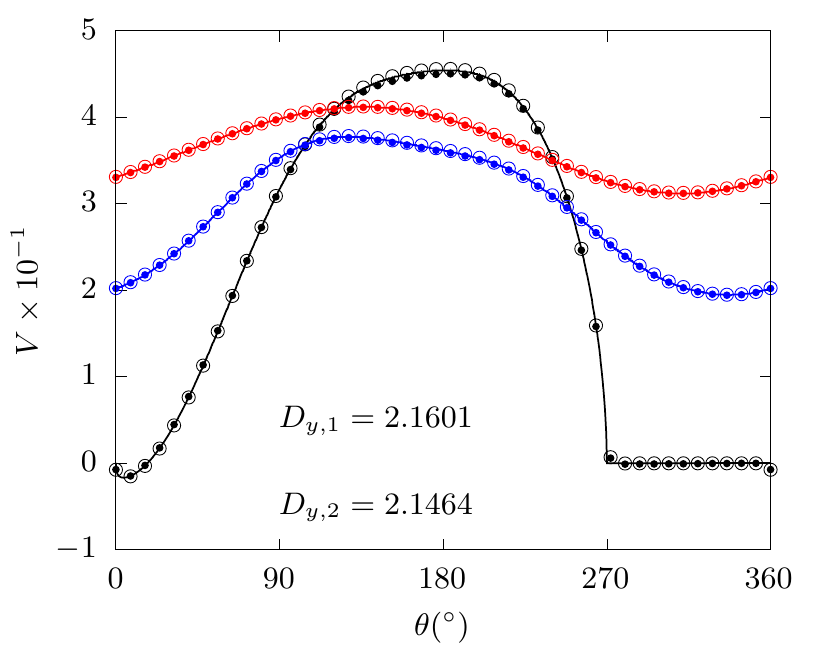}
    \\
    (e) Horizontal displacement
    &
      (f) Vertical displacement
  \end{tabular}
  \caption{Schematic diagram and stress and displacement comparisons of Case A}
  \label{fig:2}
\end{figure}

\begin{figure}[htb]
  \centering
  \begin{tabular}{cc}
      \begin{tikzpicture}
        \filldraw [gray!30] (0,0) circle [radius = 3];
        \draw [line width = 1pt, cyan] (0,0) circle [radius = 3];
        \fill [pattern = north west lines] ({3*cos(90)},{3*sin(90)}) arc [start angle = 90, end angle = -90, radius = 3] -- ({3.2*cos(-90)},{3.2*sin(-90)}) arc [start angle = -90, end angle = 90, radius = 3.2] -- ({3*cos(90)},{3*sin(90)});
        \draw [line width = 2.5pt] ({3*cos(90)},{3*sin(90)}) arc [start angle = 90, end angle = -90, radius = 3];
        \filldraw [white] (0,0) circle [radius = 0.9];
        \draw [line width = 1pt, violet] (0,0) circle [radius = 0.9];
        \foreach \x in {0,1,2,...,9} \draw [->] ({0.9*cos(360/10*\x)},{0.9*sin(360/10*\x)}) -- ({0.9*cos(360/10*\x)+0.3},{0.9*sin(360/10*\x)});
        \fill [red] ({3*cos(90)},{3*sin(90)}) circle [radius = 0.06];
        \fill [red] ({3*cos(-90)},{3*sin(-90)}) circle [radius = 0.06];
        \node at ({3*cos(-90)},{3*sin(-90)}) [red, above] {$ \theta_{1} = -\frac{\pi}{2} $};
        \node at ({3*cos(90)},{3*sin(90)}) [red, below] {$ \theta_{2} = \frac{\pi}{2} $};
        \node at (0,-0.9) [below, align=center] {$ f(\theta) = \cos\theta - {\rm i}\sin\theta $};
        \node at (0,-3.2) [below] {$ r_{i}/r_{o} $ = 0.3, $ r_{\rho}/r_{o} $ = 0.5};
      \end{tikzpicture}
    &
      \includegraphics{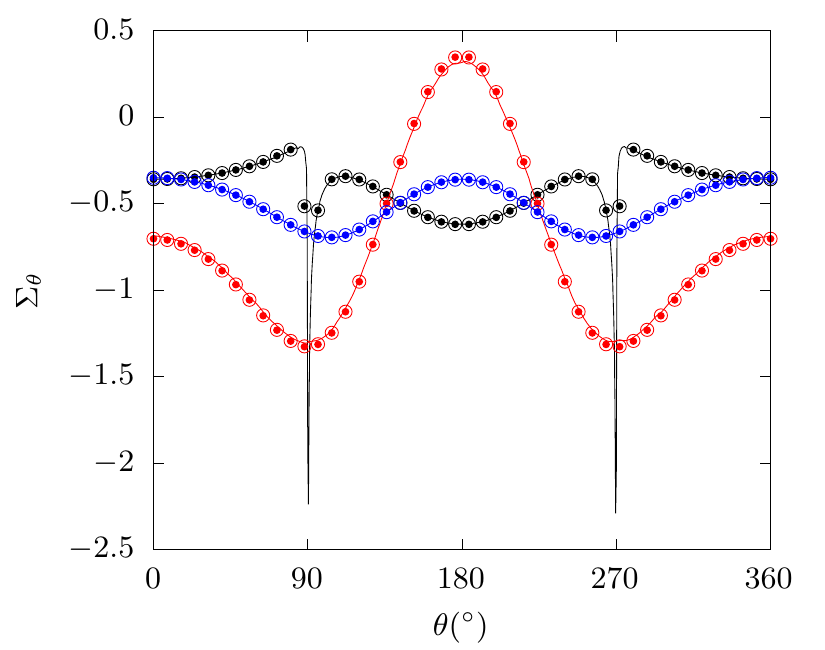}
    \\
    (a) Schematic diagram of Case B
    &
      (b) Hoop stress
    \\
    \includegraphics{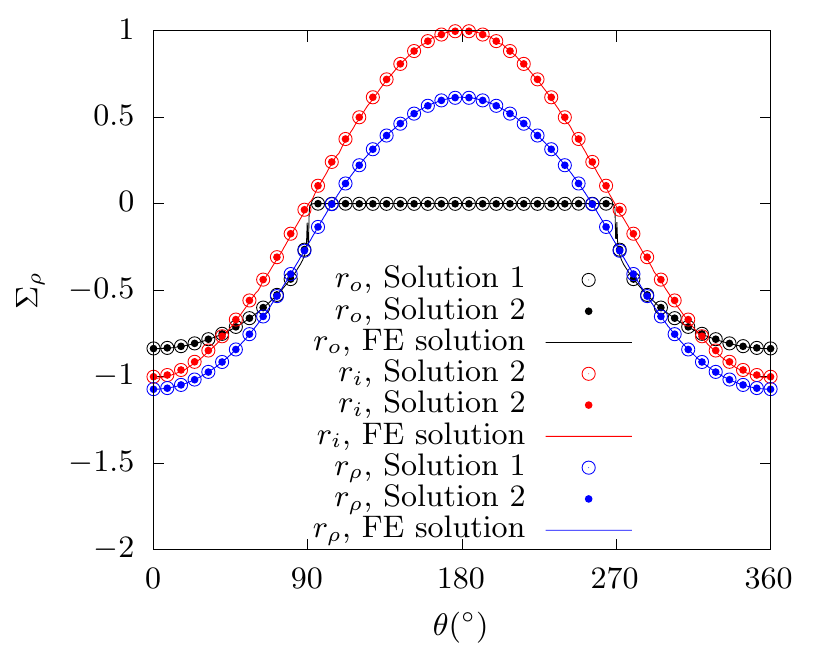}
    &
      \includegraphics{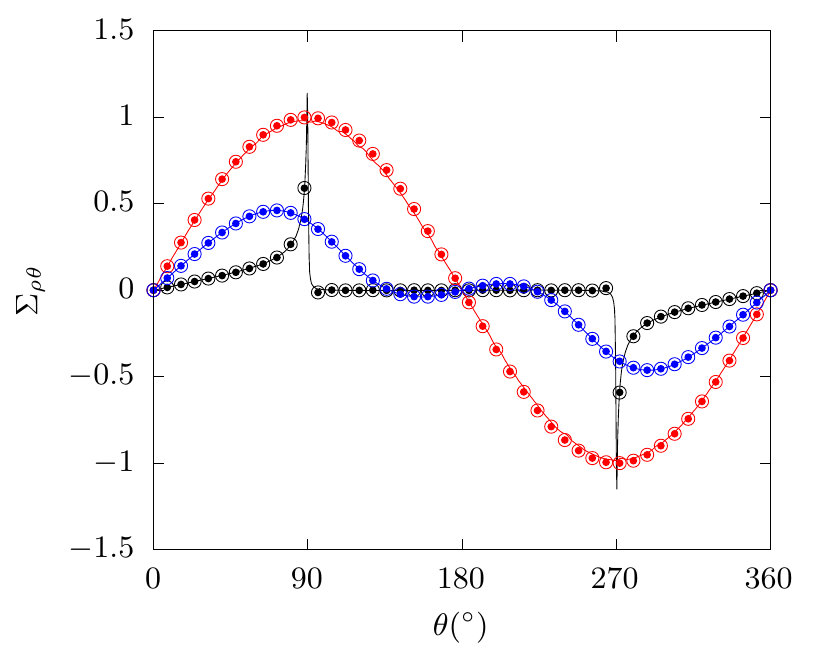}
    \\
    (c) Radial stress
    &
      (d) Shear stress
    \\
    \includegraphics{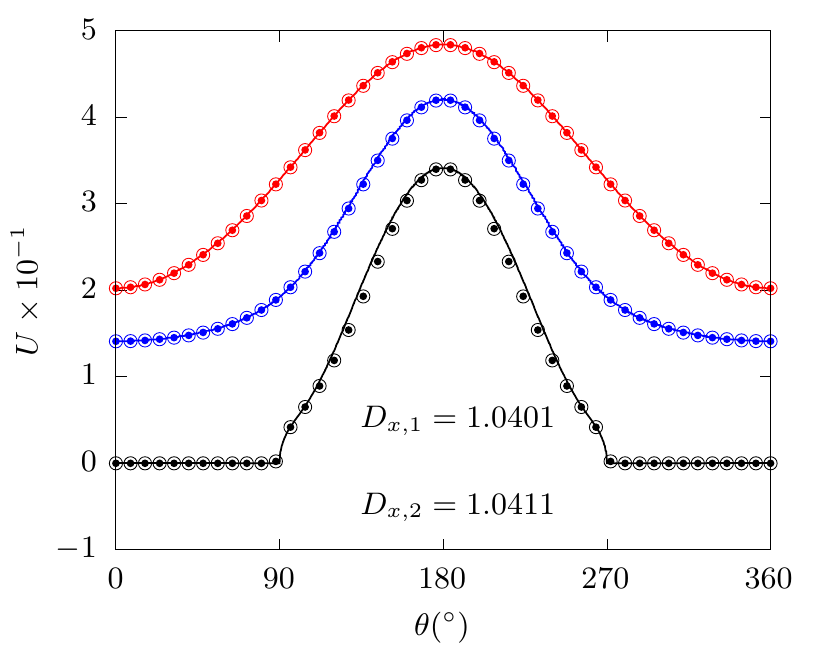}
    &
      \includegraphics{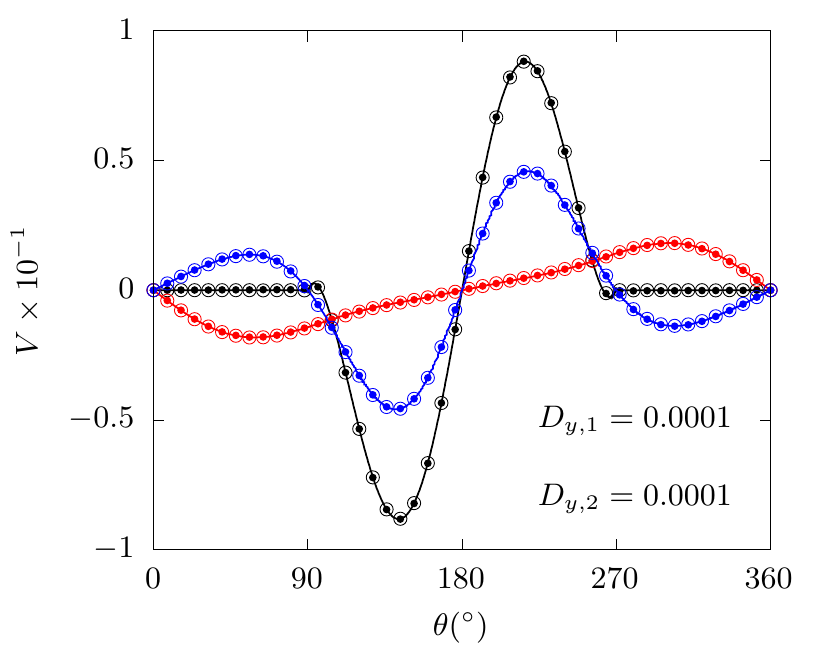}
    \\
    (e) Horizontal displacement
    &
      (f) Vertical displacement
  \end{tabular}
  \caption{Schematic diagram and stress and displacement comparisons of Case B}
  \label{fig:3}
\end{figure}

\begin{figure}[htb]
  \centering
  \begin{tabular}{cc}
      \begin{tikzpicture}
        \filldraw [gray!30] (0,0) circle [radius = 3];
        \draw [line width = 1pt, cyan] (0,0) circle [radius = 3];
        \fill [pattern = north west lines] ({3*cos(0)},{3*sin(0)}) arc [start angle = 0, end angle = -90, radius = 3] -- ({3.2*cos(-90)},{3.2*sin(-90)}) arc [start angle = -90, end angle = 0, radius = 3.2] -- ({3*cos(0)},{3*sin(0)});
        \draw [line width = 2.5pt] ({3*cos(0)},{3*sin(0)}) arc [start angle = 0, end angle = -90, radius = 3];
        \filldraw [white] (0,0) circle [radius = 1.5];
        \draw [line width = 1pt, violet] (0,0) circle [radius = 1.5];
        \foreach \x in {0,1,2,...,9} \draw [->] ({1.2*cos(360/10*\x)},{1.2*sin(360/10*\x)}) -- ({1.5*cos(360/10*\x)},{1.5*sin(360/10*\x)});
        \fill [red] ({3*cos(0)},{3*sin(0)}) circle [radius = 0.06];
        \fill [red] ({3*cos(-90)},{3*sin(-90)}) circle [radius = 0.06];
        \node at ({3*cos(-90)},{3*sin(-90)}) [red, above] {$ \theta_{1} = -\frac{\pi}{2} $};
        \node at ({3*cos(0)},{3*sin(0)}) [red, left] {$ \theta_{2} = 0 $};
        \node at (0,-1.5) [below, align=center] {$ f(\theta) = 1 $};
        \node at (0,-3.2) [below] {$ r_{i}/r_{o} $ = 0.5, $ r_{\rho}/r_{o} $ = 0.7};
      \end{tikzpicture}
    &
      \includegraphics{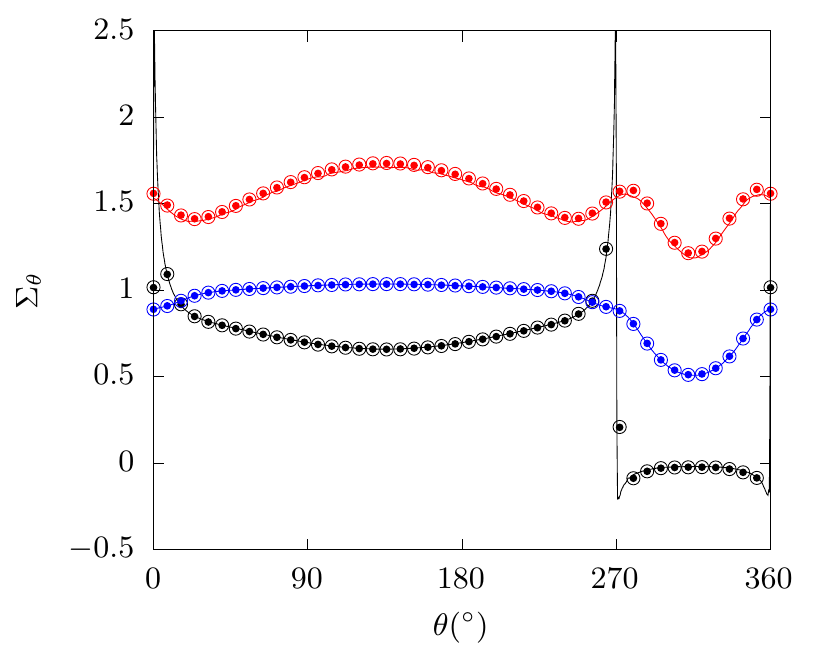}
    \\
    (a) Schematic diagram of Case C
    &
      (b) Hoop stress
    \\
    \includegraphics{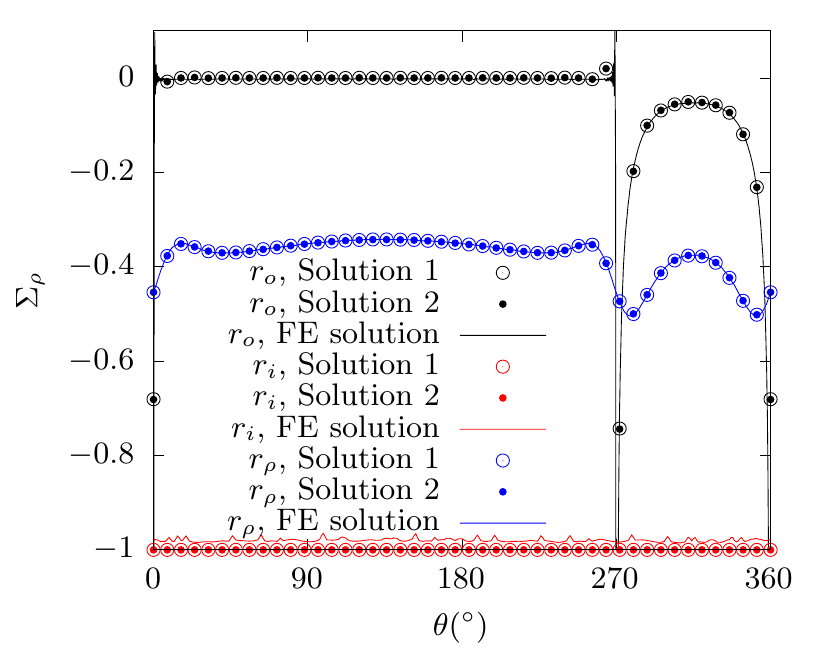}
    &
      \includegraphics{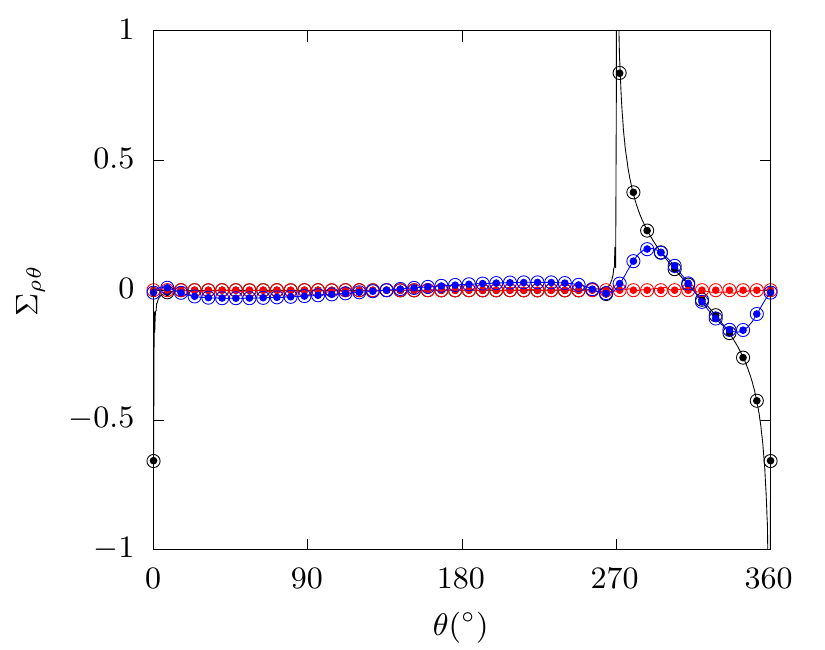}
    \\
    (c) Radial stress
    &
      (d) Shear stress
    \\
    \includegraphics{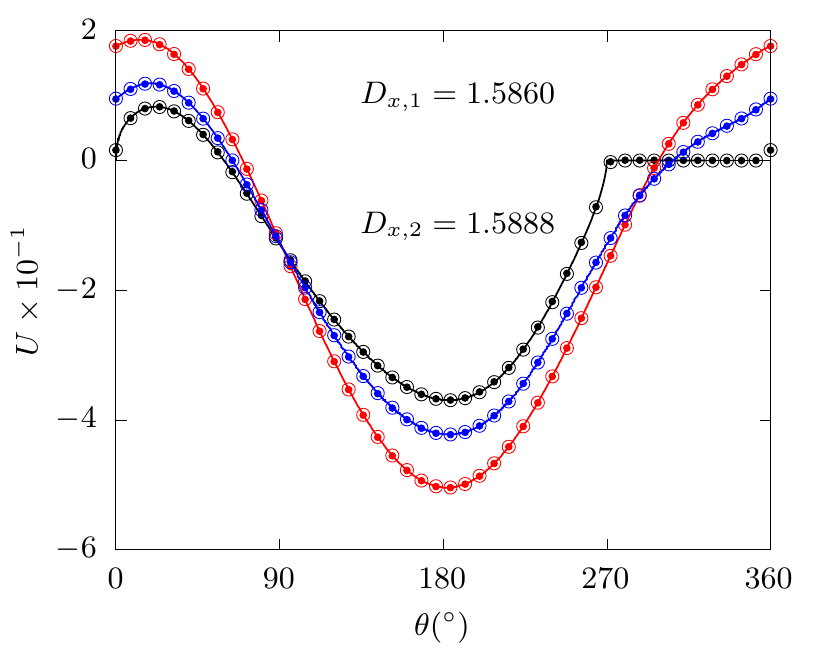}
    &
      \includegraphics{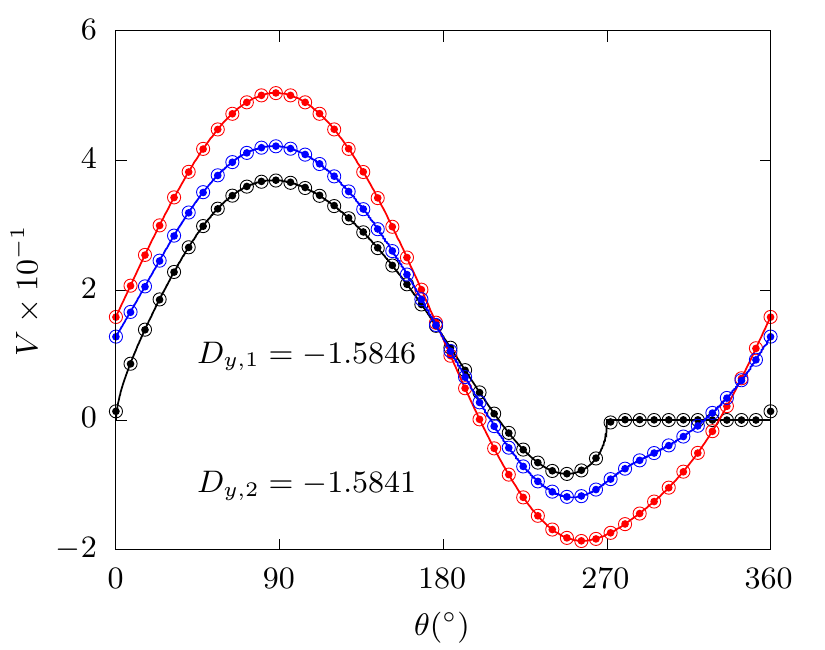}
    \\
    (e) Horizontal displacement
    &
      (f) Vertical displacement
  \end{tabular}
  \caption{Schematic diagram and stress and displacement comparisons of Case C}
  \label{fig:4}
\end{figure}

\begin{figure}[htb]
  \centering
  \begin{tabular}{cc}
      \begin{tikzpicture}
        \filldraw [gray!30] (0,0) circle [radius = 3];
        \draw [line width = 1pt, cyan] (0,0) circle [radius = 3];
        \fill [pattern = north west lines] ({3*cos(90)},{3*sin(90)}) arc [start angle = 90, end angle = -90, radius = 3] -- ({3.2*cos(-90)},{3.2*sin(-90)}) arc [start angle = -90, end angle = 90, radius = 3.2] -- ({3*cos(90)},{3*sin(90)});
        \draw [line width = 2.5pt] ({3*cos(90)},{3*sin(90)}) arc [start angle = 90, end angle = -90, radius = 3];
        \filldraw [white] (0,0) circle [radius = 2.1];
        \draw [line width = 1pt, violet] (0,0) circle [radius = 2.1];
        \foreach \x in {0,1,2,...,9} \draw [->] ({2*cos(360/10*\x)},{2*sin(360/10*\x)}) arc [start angle = 360/10*\x, end angle = 360/10*\x+30, radius = 2] ;
        \fill [red] ({3*cos(90)},{3*sin(90)}) circle [radius = 0.06];
        \fill [red] ({3*cos(-90)},{3*sin(-90)}) circle [radius = 0.06];
        \node at ({3*cos(-90)},{3*sin(-90)}) [red, above] {$ \theta_{1} = -\frac{\pi}{2} $};
        \node at ({3*cos(90)},{3*sin(90)}) [red, below] {$ \theta_{2} = \frac{\pi}{2} $};
        \node at (0,-1.8) [above, align=center] {$ f(\theta) = {\rm i} $};
        \node at (0,-3.2) [below] {$ r_{i}/r_{o} $ = 0.7, $ r_{\rho}/r_{o} $ = 0.9};
      \end{tikzpicture}
    &
      \includegraphics{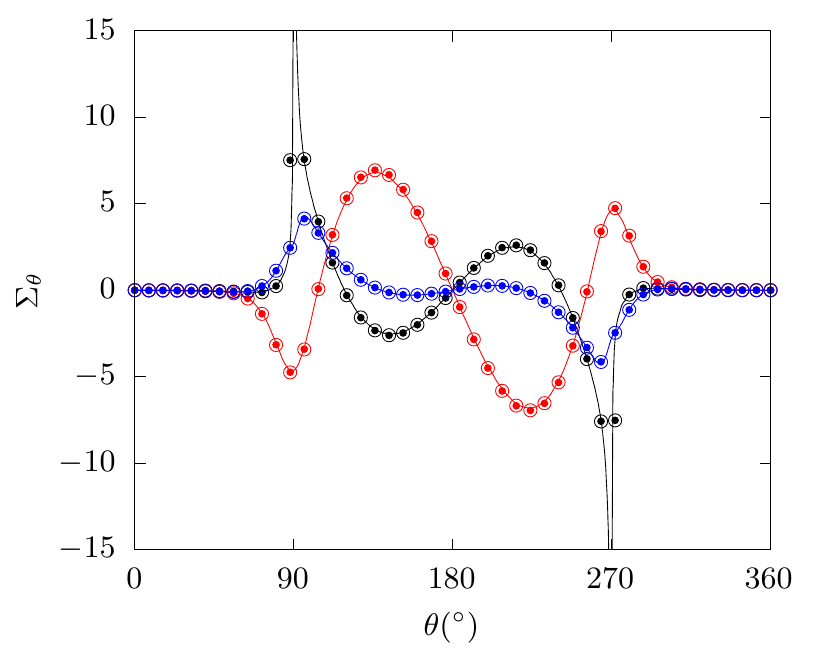}
    \\
    (a) Schematic diagram of Case D
    &
      (b) Hoop stress
    \\
    \includegraphics{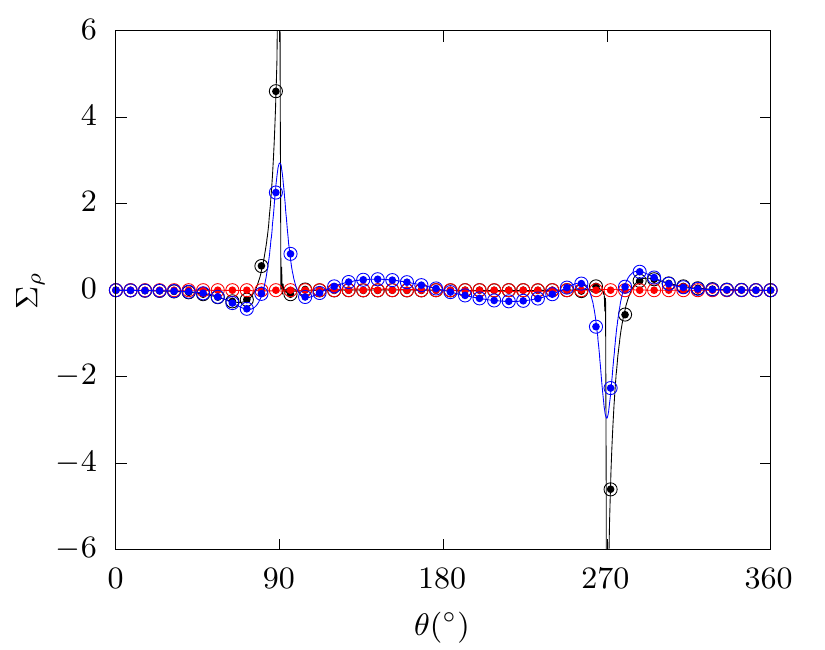}
    &
      \includegraphics{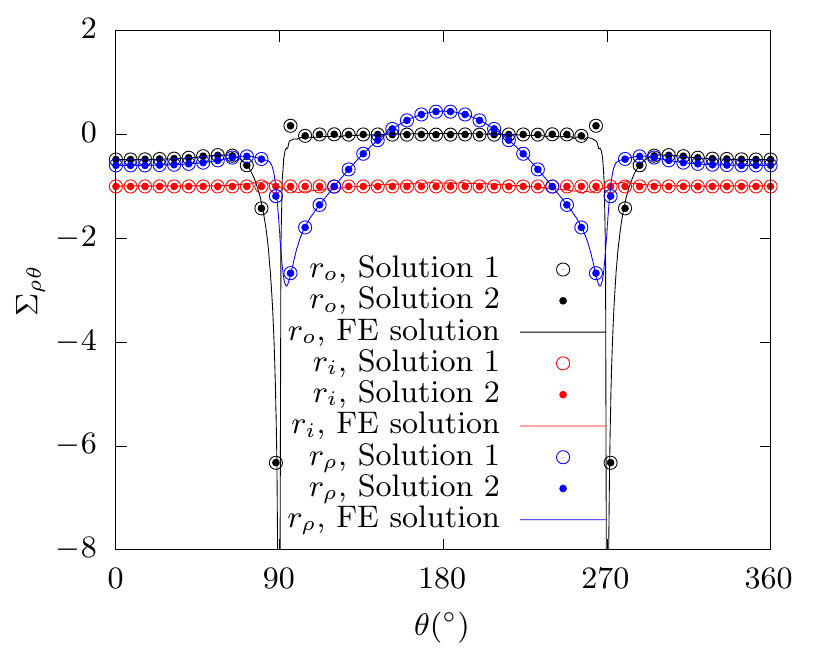}
    \\
    (c) Radial stress
    &
      (d) Shear stress
    \\
    \includegraphics{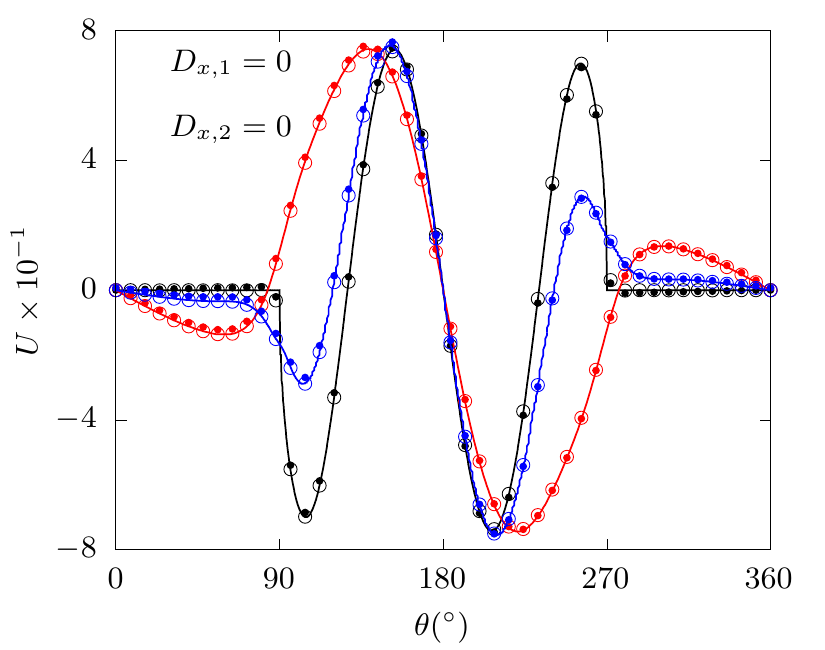}
    &
      \includegraphics{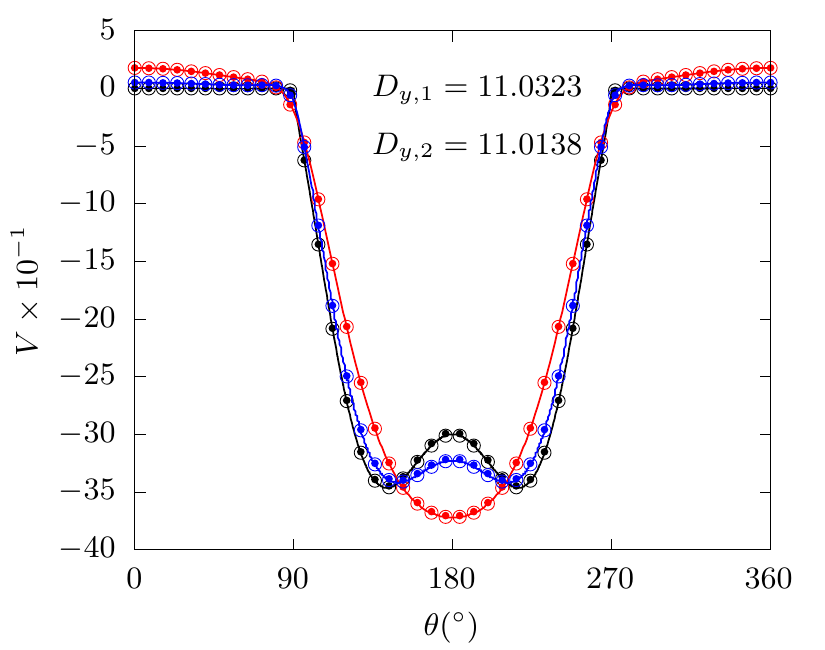}
    \\
    (e) Horizontal displacement
    &
      (f) Vertical displacement
  \end{tabular}
  \caption{Schematic diagram and stress and displacement comparisons of Case D}
  \label{fig:5}
\end{figure}

The schematic diagrams of these four cases are shown in Figs. \ref{fig:2}a-\ref{fig:5}a, respectively. Case A denotes a general case with the combination of an unbalanced traction along the inner boundary and an arbitrary support along the outer boundry. Case B is a particular case of Case A with axisymmetrical geometry, support, and traction. Case C is a general case with the combination of a balanced traction along the innner boundary and an arbitrary support along the outer boundary. Case D is a particular case of Case A with axisymmetrical geometry and support and centrosymmtric traction.

The input paramters ( $ \nu $, $ p_{k} + {\rm i} q_{k} $, $ r_{i}/r_{o} $, $ \theta_{1} $, $ \theta_2 $, $ N $), condition numbers ($ N_{C1} $, $ N_{C2} $), the maximum iteration reps ($ Q $), and $ M_{1} $ and $ M_{2} $ are listed in Table \ref{tab:1}.The truncation number takes $ N=60 $ for accuracy. Table \ref{tab:1} indicates that the condition numbers of the coefficient matrice are small, thus, the computation results are accurate.

Table \ref{tab:1} indicates that iteration reps generally increase with the ratio of the inner radius to the outer radius of the annuli. In a more detailed sense, all coefficients in Eqs. (\ref{eq:4.18}) and (\ref{eq:4.33}) can be written into sum of $ r^{k} (k \geq 2) $, except for $ (1-r^{-2})r^{2} $ in Eqs. (\ref{eq:4.18b}) and (\ref{eq:4.33b}). As $ r $ approaches 0, $ r^{k} (k \geq 2) $ would approach 0, and would converge very fast, but $ (1-r^{-2})r^{2} $ would approach $-1$, and would converge much slower; similarly, as $ r $ approaches 1, $ (1-r^{-2})r^{2} $ would approach 0, and converge fast, but $ r^{k} (k \geq 2) $ would approach 1, and would converge much slower. Such arrangements of the coefficinets reveal an insight that the successive approximation method in the solutions would be relatively slow when $ r $ approaches 0 or 1, and may reach maximum convergence speed for some value between 0 and 1. Such an insight has been verified in Table \ref{tab:1}. The iteration computation is pretty fast, and generally each case takes less than 1 sec. Thus, the computation accuracy and speed are both satisfactory. 

Substituting the same input parameters of these four cases into \texttt{ABAQUS 2020}, the finite element solutions are correspondingly obtained for comparisons. To ensure the results of the finite element solutions are accurate enough, we respectively set 360 and 720 seeds along the inner and outer peripheries of the annuli in the finite element models for all four cases, and the element quantities of these four cases are 95462, 61817, 38005, and 30007, respectively.

The comparison results between the analytical solutions and the finite element solutions for the four cases are shown in Figs. \ref{fig:2}-\ref{fig:5}, respectively. Note that the rule of signs in \texttt{ABAQUS} is different from the analytical solutions in this study. Thus, all the stress and displacement components in these four cases take a negative sign to be identical to the results obtained via corresponding finite element solutions.

Figs. \ref{fig:2}-\ref{fig:5} suggest good agreements for both stress and displacement components among Solution 1, Solution 2, and the finite element solution along all three data-selecting circles of radii $ r_{o} $, $ r_{i} $, and $ r_{\rho} $. The analytical results along the data-selecting circle of radius $ r_{i} $ in Figs. \ref{fig:2}c-\ref{fig:5}c and Figs. \ref{fig:2}d-\ref{fig:5}d respectively show complete agreements with the corresponding analytic expressions in polar form of the boundary conditions along the inner peripheries of these four annuli (a negative sign should be applied to keep the same sign rule), as is illustrated in Figs. \ref{fig:2}a-\ref{fig:5}a. The radial stress of $ r_{i} $ of the finite element result in Fig. (\ref{fig:4}c) is slightly deviated from the accurate value $-1$. The reason may be that the boundary condition of displacement constraint is strictly satisfied with prior accuracy requirement in \texttt{ABAQUS}, while the one of surface traction can be relatively relaxed with secondary accuracy requirement to ensure convergence. The comparisons above suggest Solutions 1 and 2 in this study would provide almost the same stress and displacement results, which are more accurate and robust than the finite element solution.

The consistencies of the stress and displacement components solved by Solutions 1 and 2 in Figs. \ref{fig:2}-\ref{fig:5} may provide an insight that these two solutions are the same in the numerical perspective, though these two solutions respectively employ different outer boundary conditions in Eqs. (\ref{eq:4.11}) and (\ref{eq:4.24}). Now we should examine the boundary conditions along the outer boundaries $ {\bm C}_{12} $ and $ {\bm C}_{11} $ using the numerical results to further identify and verify such an insight.

\begin{figure}[htb]
  \centering
  \begin{tabular}{cc}
    \includegraphics{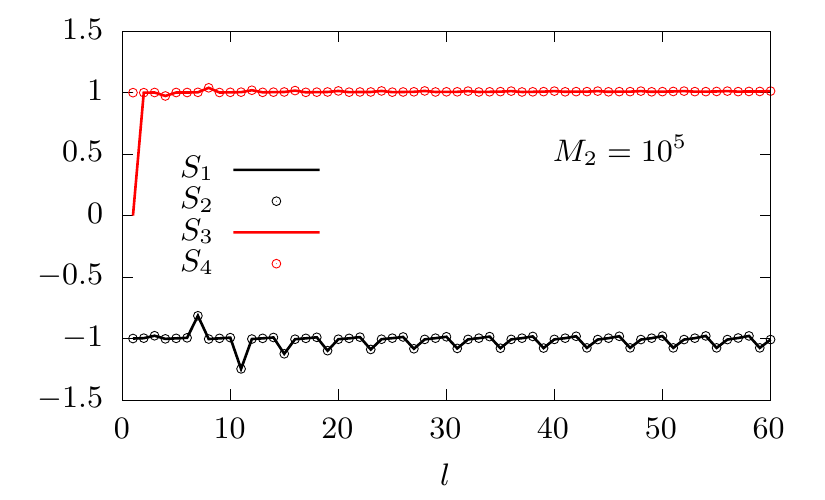}
    &
      \includegraphics{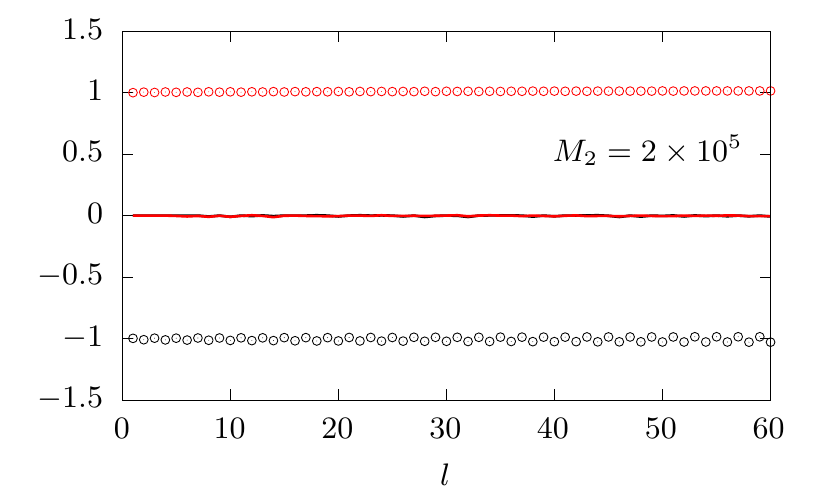}
    \\
    (a) Case A
    &
      (b) Case B
    \\
    \includegraphics{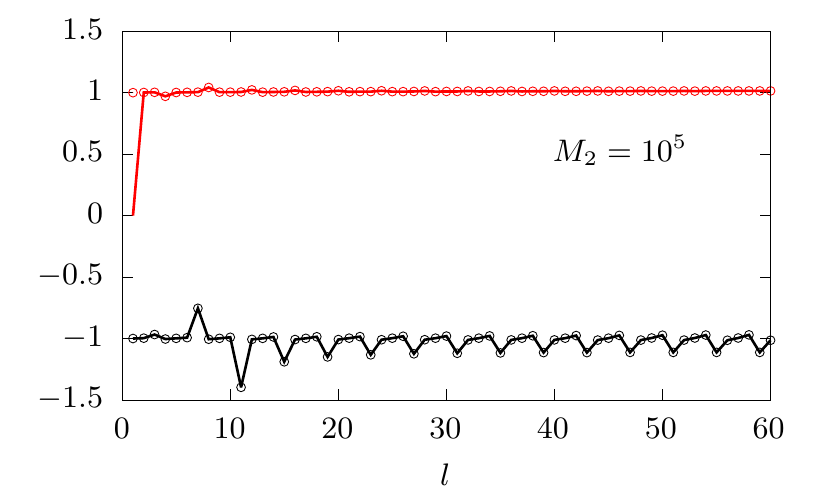}
    &
      \includegraphics{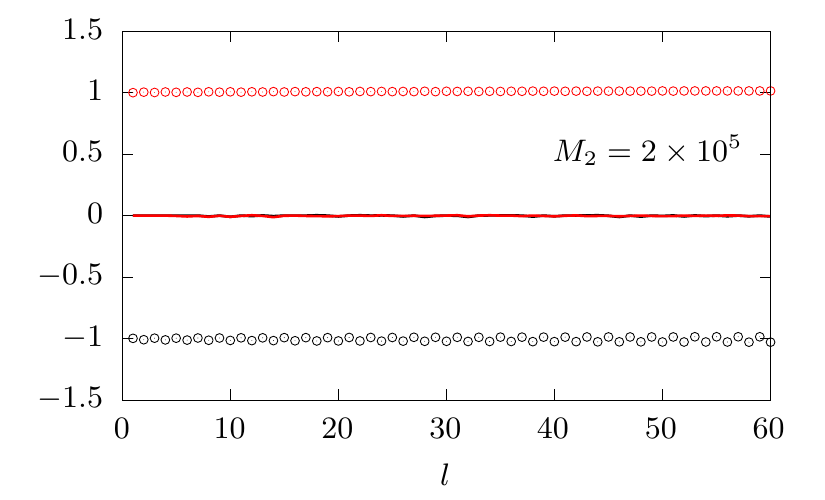}
    \\
    (c) Case C
    &
      (d) Case D
  \end{tabular}
  \caption{Comparisons of coefficients in Eq. (\ref{eq:4.27}) in Solution 2 for four cases}
  \label{fig:6}
\end{figure}

First, we examine the boundary conditions along boundary $ {\bm C}_{12} $ of these two solutions by verifying the equilities proposed in Eq. (\ref{eq:4.27}). To facilitate description, the following four variables are set:
\begin{equation*}
  S_{1} = \frac{\Re \left[ \frac{2\pi{\rm i}\kappa}{1+\kappa} \alpha_{l}\right]}{\Re [c_{12,-1-l}]}, \quad S_{2} = \frac{\Im \left[ \frac{2\pi{\rm i}\kappa}{1+\kappa} \alpha_{l}\right]}{\Im [ c_{12,-1-l}]}, \quad S_{3} = \frac{\Re \left[ \frac{2\pi{\rm i}\kappa}{1+\kappa} \beta_{l}\right]}{\Re [ c_{12,-1+l}]}, \quad S_{4} = \frac{\Im \left[ \frac{2\pi{\rm i}\kappa}{1+\kappa} \beta_{l} \right]}{\Im [ c_{12,-1+l} ]}
\end{equation*}
According to Eq. (\ref{eq:4.27}), $ S_{1} $ and $ S_{2} $ should be theoretically equal to $ -1 $, while $ S_{3} $ and $ S_{4} $ should be theoretically equal to 1. Since $ c_{12,k} $ are computed via numerical integrals in Eq. (\ref{eq:B.4}), and would be sensitive to accuracy, thus, we magnify the value of $ M_{2} $ in Table \ref{tab:1} by ten. The computation results of these four variables in four cases are presented in Fig. \ref{fig:6}. Apparently, the results in Fig. \ref{fig:6} are approximate to the expected values. The zero values indicate symmetry. If a larger value of $ M_{2} $ is used, the oscillation would be gradually eliminated. Therefore, the numerical results in Fig. \ref{fig:6} validates Eq. (\ref{eq:4.27}), and subsequently determine the consistency of Eqs. (\ref{eq:4.7a}) and (\ref{eq:4.23}).

\begin{figure}[htb]
  \centering
  \begin{tabular}{cc}
    \includegraphics{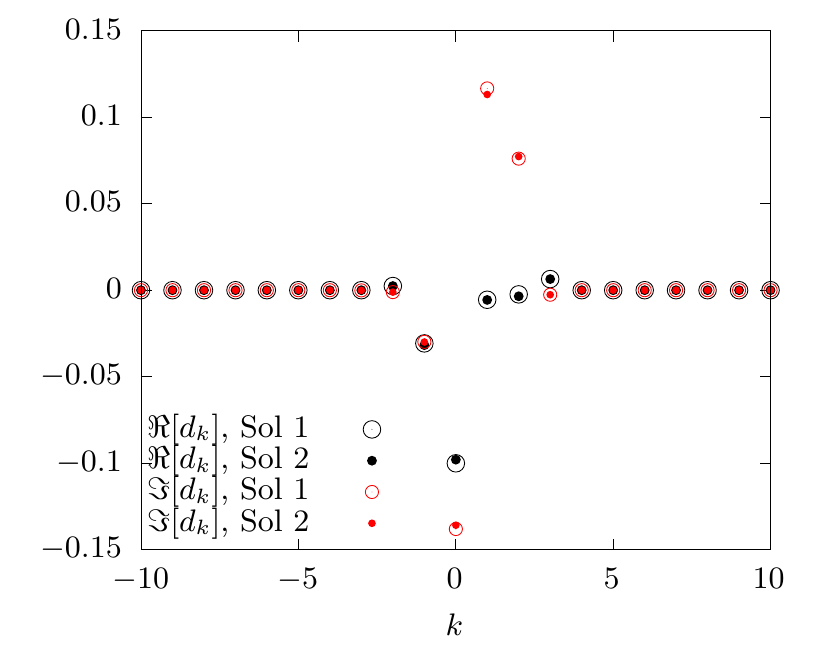}
    &
      \includegraphics{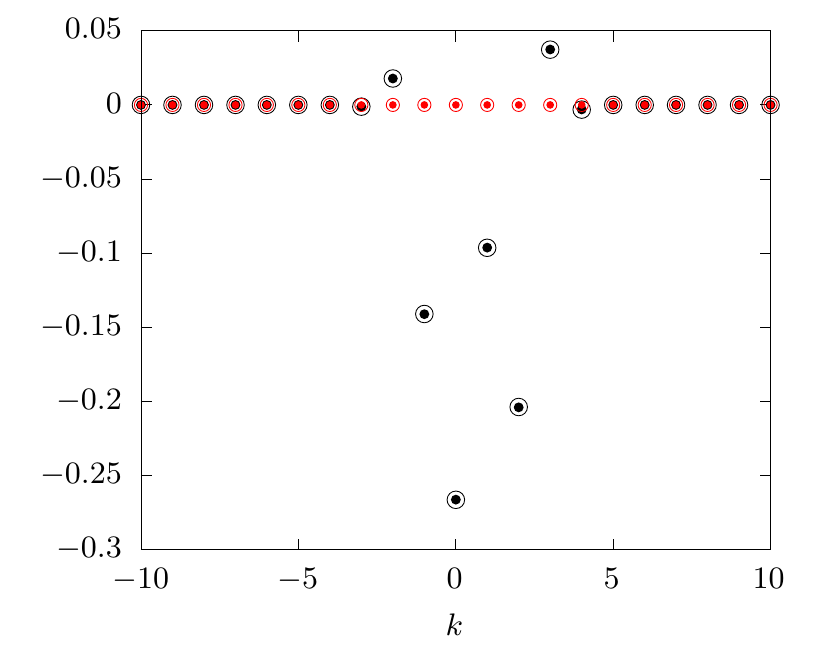}
    \\
    (a) Case A
    &
      (b) Case B
    \\
    \includegraphics{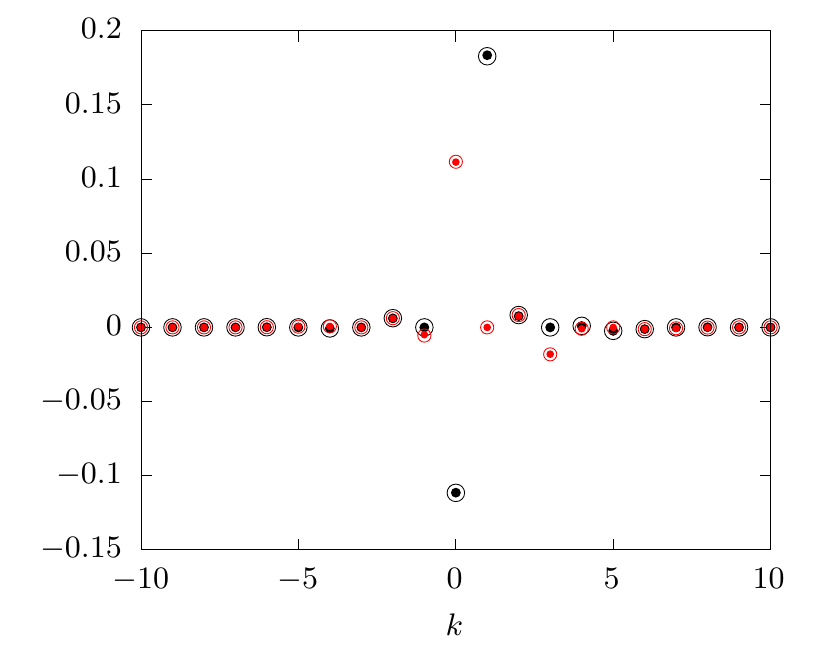}
    &
      \includegraphics{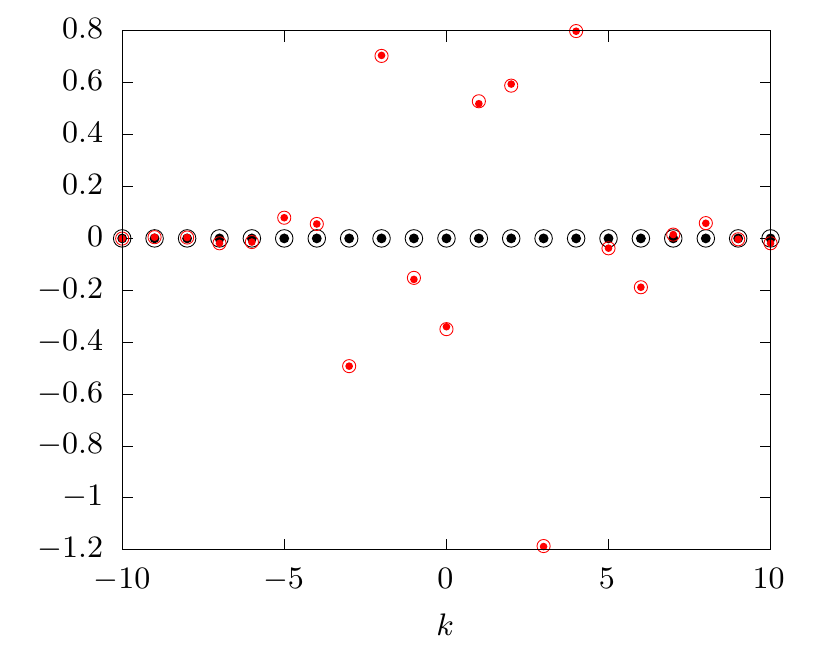}
    \\
    (c) Case C
    &
      (d) Case D
  \end{tabular}
  \caption{Comparisons of real and imaginary parts of $ d_{k} $ between Solutions 1 and 2 for four cases}
  \label{fig:7}
\end{figure}

Further, we examine the boundary conditions along boundary $ {\bm C}_{11} $ by comparing the values of $ d_{k} $ solved by Solutions 1 and 2. Solutions 1 and 2 both employ boundary conditions in Eq. (\ref{eq:4.7}), while respectively employ Eqs. (\ref{eq:4.11}) and (\ref{eq:4.24}). If the values of $ d_{k} $ solved by these two solutions are the same, the boundary conditions in Eqs. (\ref{eq:4.11}) and (\ref{eq:4.24}) would be equivalent. Fig. \ref{fig:7} shows the comparisons of the real and imaginary parts of $ d_{k} $ solved by Solutions 1 and 2, and the results suggest highly identities of $ d_{k} $ solved by these two solutions, indicating that Eq. (\ref{eq:4.24}) is a linear combination of Eqs. (\ref{eq:4.7}) and (\ref{eq:4.11}) from a linear algebra perspective, and that these two solutions are numerically equivalent.

\section{Remarks and discussion}
\label{sec:remark-discu}

(a) The definition domains of $ \varphi^{\prime+}(z) $ and $ \varphi^{\prime-}(z) $ are breached in Solution 1 to apply the residue theorem, and a mathematically elegant and simple solution method is established. To be strict, the breaching of definition domains may be mathematically flawed. Meanwhile, the definition domains of $ \varphi^{\prime+}(z) $ and $ \varphi^{\prime-}(z) $ remain intact in Solution 2. Such a treatment shows more mathematical strictness coupled with a more complicated solution method. Comparing to Solution 1, each iteration rep in Solution 2 is divided into two subiterations to ensure convergence. Furthermore, the improper integrals in Eq. (\ref{eq:4.25}) is numerically computed. Thus, Solution 2 would require more computation intensity than Solution 1, and the coding labor of Solution 2 is also much more than that of Solution 1 in the numercal cases. The four numerical cases verify in detail that these two parallel solutions show mutually numerical equivalence. Therefore, Solution 1 may be preferentially considered to lower computation intensity and coding labor for practical use, if mathematical strictness is not in priority.

(b) The deduction in Eqs. (\ref{eq:4.9})-(\ref{eq:4.11}) eliminates the ambiguities of the similar procedures in Ref \cite[]{Sugiura1973TheMB, sugiura1969heavy}. The unknown coefficients $ d_{k} $ in Eq. (\ref{eq:4.1}) should and should only be determined according to the boundary conditions, instead of other conditions, because such a procedure is the correct and strict procedure between a general solution and a determined solution. It seems that Eq. (23) of Ref \cite[]{sugiura1969heavy} and Eq. (33) of Ref \cite[]{Sugiura1973TheMB} are derived from the genreal relationship between displacements and complex potentials. To be more specific and convenient, we illustrate the procedure using the symbols in this paper. Substituting Eqs. (\ref{eq:4.4a}) and (\ref{eq:4.5}) into Eq. (\ref{eq:2.1}) yields:
\begin{equation*}
  \begin{aligned}
    g(z) = & \kappa \int \sum\limits_{k=-\infty}^{\infty} A_{k}z^{k} {\rm d}z - z \sum\limits_{k=-\infty}^{\infty} \overline{A}_{k}\overline{z}^{k} - \int \sum\limits_{k=-\infty}^{\infty} \left[ B_{-k-2} - (k+1) \overline{A}_{k+2} \right] \overline{z}^{k} {\rm d}\overline{z} \\
    = & \kappa A_{-1} {\rm Log}z - B_{-1} {\rm Log}z + g_{0}(z,\overline{z})  \\
    = & (\kappa A_{-1} + B_{-1})\cdot 2\pi{\rm i} n + \kappa A_{-1} \log z - B_{-1} \log z + g_{0}(z,\overline{z})
  \end{aligned}
\end{equation*}
where $ {\rm Log} $ and $ \log $ denote the multi-value and single-value natural logarithmic functions, respectively, $ g_{0}(z,\overline{z}) $ denotes the rest single-value functions. Apparently, the single-valueness of displacement in the annulus requires $\kappa A_{-1} + B_{-1} = 0$, which coincides with Eq. (\ref{eq:4.11}) consequently. So it seems that the deduction above is correct. Whereas it is wrong in conception level. Since $ A_{-1} $ and $ B_{-1} $ have been denoted as intermediate symbols to obtain the solution of $ d_{k} $ in Eq. (\ref{eq:4.4}), they belong to the procedure between the general solution and the determined solution in Eq. (\ref{eq:4.1})-(\ref{eq:4.23}), indicating that only the boundary conditions should be used. Whereas, the deduction above is not related to any boundary condition. Thus, the deduction above can not be used in the solving procedure.

(c) Comparing to the solutions in Ref \cite[]{Sugiura1973TheMB, sugiura1969heavy, yau1968mixed}, the mixed boundary conditions in this paper consist of a partially fixed constraint acting along the outer periphery and an arbitrary traction acting along the inner periphery. Thus, the parallel solutions in this paper are both generalized and can be potentially used in more complicated and real situations.

(d) Both solutions in this paper have been verified via numerical cases by comparing to corresponding finite element results. Furthermore, the results reveal that the proposed solution has much higher accuracy and better robustness than the corresponding finite element solution.

\appendix
\clearpage
\section{Appendix A}
\label{sec:A}

The analytic expressions of the coefficients in the Taylor's expansions are expressed as:
\begin{equation}
  \label{eq:A.1}
  \begin{cases}
    \alpha_{0} = & - t_{1}^{-\gamma} t_{2}^{\gamma-1} \\
    \alpha_{1} = & - t_{1}^{-\gamma} t_{2}^{\gamma-1} [\gamma t_{1}^{-1} - (\gamma-1) t_{2}^{-1}] \\
    \alpha_{k} = & - t_{1}^{-\gamma} t_{2}^{\gamma-1}[ (-1)^{k} \frac{h(-\gamma,k)}{k!} t_{1}^{-k} + (-1)^{k} \frac{h(\gamma-1,k)}{k!} t_{2}^{-k} \\
                 & + (-1)^{k} \sum\limits_{l=1}^{k-1} \frac {h(-\gamma,l)} {l!} \frac {h(\gamma-1,k-l)} {(k-l)!} \cdot t_{1}^{-l} t_{2}^{-k+l}], \quad k \geq 2
  \end{cases}
\end{equation}

\begin{equation}
  \label{eq:A.2}
  \begin{cases}
    \beta_{1} = & 1 \\
    \beta_{2} = & \gamma t_{1} - (\gamma-1) t_{2} \\
    \beta_{k} = & (-1)^{k-1} \frac{h(-\gamma,k-1)}{(k-1)!} t_{1}^{k-1} + (-1)^{k-1} \frac{h(\gamma-1,k-1)} {(k-1)!} t_{2}^{k-1} \\
                & + (-1)^{k-1} \sum\limits_{l=1}^{k-2} \frac {h(-\gamma,l)} {l!} \frac {h(\gamma-1,k-l-1)} {(k-l-1)!} \cdot t_{1}^{l} t_{2}^{k-l-1}, \quad k \geq 3
  \end{cases}
\end{equation}
where
\begin{equation}
  \label{eq:A.3}
  \begin{cases}
    h(-\gamma,k) = -\gamma (-\gamma-1) (-\gamma-2) \cdots (-\gamma-k+1) \\
    h(\gamma-1,k) = (\gamma-1) (\gamma-2) (\gamma-3) \cdots (\gamma-k) \\
  \end{cases}
\end{equation}

\clearpage
\section{Appendix B}
\label{sec:B}

To avoid the multi-valuedness of the improper integrals in Eq. (\ref{eq:4.26}), we would conduct the integrals in real domain. Thus, the integrands in Eq. (\ref{eq:4.26}) should be prepared into the following form:
\begin{equation}
  \label{eq:B.1}
  \begin{aligned}
    & ({\rm e}^{{\rm i}\theta}-{\rm e}^{{\rm i}\theta_{1}})^{-\gamma} ({\rm e}^{{\rm i}\theta}-{\rm e}^{{\rm i}\theta_{2}})^{\gamma-1} \cdot {\rm e}^{{\rm i}(k+1)\theta} \\
    = & K_{0} \left( \sin\frac{\theta-\theta_{1}}{2} \cdot \sin\frac{\theta-\theta_{2}}{2} \right)^{-\frac{1}{2}} \left( \frac{\sin\frac{\theta-\theta_{2}}{2}}{\sin\frac{\theta-\theta_{1}}{2}} \right)^{{\rm i}\lambda} \cdot {\rm e}^{{\rm i}\left(k+\frac{1}{2}\right)\theta} \\
    =
    & \left\{
    \begin{aligned}
       - {\rm i} {\rm e}^{\pi\lambda} & K_{0} \left( \sin\frac{\theta-\theta_{1}}{2} \cdot \sin\frac{\theta_{2}-\theta}{2} \right)^{-\frac{1}{2}} \cdot \left[ \cos \eta_{2}(k,\theta) + {\rm i} \sin \eta_{2}(k,\theta) \right], \quad \theta \in (\theta_{1}, \theta_{2}) \\
      & K_{0} \left( \sin\frac{\theta-\theta_{1}}{2} \cdot \sin\frac{\theta-\theta_{2}}{2} \right)^{-\frac{1}{2}} \cdot \left[ \cos \eta_{1}(k,\theta) + {\rm i} \sin \eta_{1}(k,\theta) \right], \quad \theta \in (\theta_{2}, \theta_{1}+2\pi) \\
    \end{aligned}
    \right.
  \end{aligned}
\end{equation}
where
\begin{equation}
  \label{eq:B.2}
  \left\{
    \begin{aligned}
      K_{0} & = \displaystyle - \frac{{\rm i}}{2} \exp\left[\frac{\lambda}{2}(\theta_{1}-\theta_{2})-\frac{{\rm i}}{4}(\theta_{1}+\theta_{2})\right] \\
      \eta_{2}(k,\theta) & = \displaystyle \left(k+\frac{1}{2}\right)\theta + \lambda \ln \frac{\sin\frac{\theta_{2}-\theta}{2}} {\sin\frac{\theta-\theta_{1}}{2}} \\
      \eta_{1}(k,\theta) & = \displaystyle \left(k+\frac{1}{2}\right)\theta + \lambda \ln \frac{\sin\frac{\theta-\theta_{2}}{2}} {\sin\frac{\theta-\theta_{1}}{2}} \\
    \end{aligned}
  \right.
\end{equation}
Then the improper interals in Eq. (\ref{eq:4.26}) can be equivalently modified as
\begin{subequations}
  \label{eq:B.3}
  \begin{equation}
    \label{eq:B.3a}
    \begin{aligned}
      c_{12,k} = & - {\rm e}^{\pi\lambda} K_{0} \lim_{\delta \rightarrow 0} \int_{\theta_{1}+\delta}^{\theta_{2}-\delta} \left( \sin\frac{\theta-\theta_{1}}{2} \cdot \sin\frac{\theta_{2}-\theta}{2} \right)^{-\frac{1}{2}} \cdot \cos\eta_{2}(k,\theta) {\rm d}\theta \\
                 & - {\rm i} {\rm e}^{\pi\lambda} K_{0} \lim_{\delta \rightarrow 0} \int_{\theta_{1}+\delta}^{\theta_{2}-\delta} \left( \sin\frac{\theta-\theta_{1}}{2} \cdot \sin\frac{\theta_{2}-\theta}{2} \right)^{-\frac{1}{2}} \sin\eta_{2}(k,\theta) {\rm d}\theta
    \end{aligned}
  \end{equation}
  \begin{equation}
    \label{eq:B.3b}
    \begin{aligned}
      c_{11,k} = & {\rm i} K_{0} \lim_{\delta \rightarrow 0} \int_{\theta_{2}+\delta}^{\theta_{1}+2\pi-\delta} \left( \sin\frac{\theta-\theta_{1}}{2} \cdot \sin\frac{\theta-\theta_{2}}{2} \right)^{-\frac{1}{2}} \cdot \cos \eta_{1}(k,\theta) {\rm d}\theta \\
                 & - K_{0} \lim_{\delta \rightarrow 0} \int_{\theta_{2}+\delta}^{\theta_{1}+2\pi-\delta} \left( \sin\frac{\theta-\theta_{1}}{2} \cdot \sin\frac{\theta-\theta_{2}}{2} \right)^{-\frac{1}{2}} \cdot \sin \eta_{1}(k,\theta) {\rm d}\theta
    \end{aligned}
  \end{equation}
\end{subequations}
where $ \delta $ denotes a small numeric. The integration direction of $ c_{12,k} $ is reversed to facilitate computation. Then the complex integrals in Eq. (\ref{eq:B.3}) turn to two real integrals free from branches of complex variable, respectively, and can be approximately obtained as
\begin{subequations}
  \label{eq:B.4}
  \begin{equation}
    \label{eq:B.4a}
    \begin{aligned}
      & \lim_{\delta \rightarrow 0} \int_{\theta_{1}+\delta}^{\theta_{2}-\delta} \left( \sin\frac{\theta-\theta_{1}}{2} \cdot \sin\frac{\theta_{2}-\theta}{2} \right)^{-\frac{1}{2}} \cos\eta_{2}(k,\theta) {\rm d}\theta \\
      = & \sum\limits_{m=1}^{M_{2}} \left( \sin\frac{\theta_{(m)}-\theta_{1}}{2} \cdot \sin\frac{\theta_{2}-\theta_{(m)}}{2} \right)^{-\frac{1}{2}} \cdot \cos\eta_{2}(k,\theta_{(m)}) \cdot \Delta\theta_{2}
    \end{aligned}
  \end{equation}
  \begin{equation}
    \label{eq:B.4b}
    \begin{aligned}
      & \lim_{\delta \rightarrow 0} \int_{\theta_{1}+\delta}^{\theta_{2}-\delta} \left( \sin\frac{\theta-\theta_{1}}{2} \cdot \sin\frac{\theta_{2}-\theta}{2} \right)^{-\frac{1}{2}} \cdot \sin\eta_{2}(k,\theta_{(m)}) {\rm d}\theta \\
      = & \sum\limits_{m=1}^{M_{2}} \left( \sin\frac{\theta_{(m)}-\theta_{1}}{2} \cdot \sin\frac{\theta_{2}-\theta_{(m)}}{2} \right)^{-\frac{1}{2}} \cdot \sin\eta_{2}(k,\theta_{(m)}) \cdot \Delta\theta_{2}
    \end{aligned}
  \end{equation}
\end{subequations}
where
\begin{equation}
  \label{eq:B.5}
  \begin{cases}
    \Delta\theta_{2} = \displaystyle \frac{\theta_{2}-\theta_{1}}{M_{2}+1} \\
    \theta_{(m)} = \theta_{1} + m \cdot \Delta\theta_{2} \\
  \end{cases}
\end{equation}
\begin{subequations}
  \label{eq:B.6}
  \begin{equation}
    \label{eq:B.6a}
    \begin{aligned}
      & \lim_{\delta \rightarrow 0} \int_{\theta_{2}+\delta}^{\theta_{1}+2\pi-\delta} \left( \sin\frac{\theta-\theta_{1}}{2} \cdot \sin\frac{\theta-\theta_{2}}{2} \right)^{-\frac{1}{2}} \cdot \cos\eta_{1}(k,\theta) {\rm d}\theta_{1} \\
      = & \sum\limits_{m=1}^{M_{1}} \left( \sin\frac{\theta_{(m)}-\theta_{1}}{2} \cdot \sin\frac{\theta_{(m)}-\theta_{2}}{2} \right)^{-\frac{1}{2}} \cdot \cos\eta_{1}(k,\theta_{(m)}) \cdot \Delta\theta_{1}
    \end{aligned}
  \end{equation}
  \begin{equation}
    \label{eq:B.6b}
    \begin{aligned}
      & \lim_{\delta \rightarrow 0} \int_{\theta_{2}+\delta}^{\theta_{1}+2\pi-\delta} \left( \sin\frac{\theta-\theta_{1}}{2} \cdot \sin\frac{\theta-\theta_{2}}{2} \right)^{-\frac{1}{2}} \cdot \sin\eta_{1}(k,\theta) {\rm d}\theta \\
      = & \sum\limits_{m=1}^{M_{1}} \left( \sin\frac{\theta_{(m)}-\theta_{1}}{2} \cdot \sin\frac{\theta_{(m)}-\theta_{2}}{2} \right)^{-\frac{1}{2}} \cdot \sin\eta_{1}(k,\theta_{(m)}) \cdot \Delta\theta_{1}
    \end{aligned}
  \end{equation}
\end{subequations}
where
\begin{equation}
  \label{eq:B.7}
  \begin{cases}
    \Delta\theta_{1} = \displaystyle \frac{\theta_{1}-\theta_{2}+2\pi}{M_{1}+1} \\
    \theta_{(m)} = \theta_{2} + m \cdot \Delta\theta_{1} \\
  \end{cases}
\end{equation}
$ M_{2} $ and $ M_{1} $ are two large positive integers. Note that the sums in Eqs. (\ref{eq:B.4}) and (\ref{eq:B.6}) start from item 1 to $ M_{2} $ and $ M_{1} $, respectively, thus, the poles $ t_{1} $ and $ t_{2} $ would not be included. As long as $ M_{2} $ and $ M_{1} $ are large enough, the approximation of Eqs. (\ref{eq:B.4}) and (\ref{eq:B.6}) would be accurate enough.

\clearpage
\section*{Acknowlegement}
\label{sec:acknowlegement}

This study is financially supported by the Natural Science Foundation of Fujian Province, China (Grant No. 2022J05190), the Scientific Research Foundation of Fujian University of Technology (Grant No. GY-Z20094), the National Natural Science Foundation of China (Grant No. 52178318), and the Education Foundation of Fujian Province (Grant No. JAT210287). The authors would like to thank Professor Changjie Zheng, Ph.D. Yiqun Huang, and Associate Professor Xiaoyi Zhang for their suggestions on this study.

\clearpage
\begin{table}[htb]
  \centering
  \begin{threeparttable}
    \begin{tabular}{ccccccccccccccc}
      \toprule
      Case & $ \nu $ & $ p_{k}+{\rm i}q_{k} $ & $r_{i}/r_{o}$ & $\theta_{1}$ & $ \theta_{2} $ & $N$ & $ N_{C1,1} $ & $ N_{C2,1} $ & $Q_{1}$ \tnote{1} & $ N_{C1,2} $ & $ N_{C2,2} $ & $Q_{2}$ \tnote{2} & $ M_{1} $ & $ M_{2} $ \\
      \midrule
      A & 0.3 & $ q_{-1} = 1 $ & 0.1 & $ -\frac{\pi}{2} $ & 0 & 60 & 8.76 & 11.92 & 140 & 8.70 & 11.89 & 139 & $ 3 \cdot 10^{4} $ & $ 1 \cdot 10^{4} $ \\
      B & 0.3 & $ p_{-1} = 1 $ & 0.3 & $ -\frac{\pi}{2} $ & $ \frac{\pi}{2} $ & 60 & 7.59 & 7.63 & 29 & 9.92 & 9.92 & 30 & $ 2 \cdot 10^{4} $ & $ 2 \cdot 10^{4} $ \\
      C & 0.3 & $ p_{0} = 1 $ & 0.5 & $ -\frac{\pi}{2} $ & 0 & 60 & 8.76 & 11.92 & 79 & 9.32 & 9.32 & 239 & $ 3 \cdot 10^{4} $ & $ 1 \cdot 10^{4} $ \\
      D & 0.3 & $ q_{0} = 1 $ & 0.7 & $ -\frac{\pi}{2} $ & $ \frac{\pi}{2} $ & 60 & 7.59 & 7.63 & 286 & 7.52 & 7.59 & 286 & $ 2 \cdot 10^{4} $ & $ 2 \cdot 10^{4} $ \\
      \bottomrule
    \end{tabular}
    \begin{tablenotes}
      \footnotesize
    \item [1] $ N_{C1,1} $, $ N_{C2,1} $, and $ Q_{1} $ denote condition numbers and maximum iteration reps for Solution 1.
    \item [2] $ N_{C1,2} $, $ N_{C2,2} $, and $ Q_{2} $ denote condition numbers and maximum iteration reps for Solution 2.
    \end{tablenotes}
  \end{threeparttable}
  \caption{Paramters of four cases}
  \label{tab:1}
\end{table}

\clearpage

\clearpage

\clearpage

\clearpage

\clearpage

\clearpage

\clearpage

\clearpage
\bibliographystyle{plainnat}
\bibliography{./citation}  %%% Uncomment this line and comment out the ``thebibliography'' section below to use the external .bib file (using bibtex) .

\end{document}